\documentclass{article}

\title{\LARGE \textbf{A Size Upper Bound for Dominating Cycles}}
\author{Zh.G. Nikoghosyan}

\begin{document}

\maketitle

\begin{abstract}
Recently it was shown (by the author) that  every graph of size $q$ (the number of edges) and minimum degree $\delta$ is hamiltonian if $q\le\delta^2+\delta-1$ (arXiv:1107.2201v1). In this paper we present the exact analog of this result for dominating cycles: if $G$ is a 2-connected graph with  $q\le8$ if $\delta=2$ and $q\le (3(\delta-1)(\delta+2)-1)/2$ if $\delta\ge3$, then each longest cycle in $G$ is a dominating cycle. The result is sharp in all respects.\\


\end{abstract}

\section{Introduction}

Only finite undirected graphs without loops or multiple edges are considered. We reserve $n$, $q$, $\delta$ and $\kappa$  to denote the number of vertices (order), the number of edges (size), the minimum degree and the connectivity of a graph, respectively. A good reference for any undefined terms is \cite{[1]}. 

The earliest sufficient condition for a graph to be hamiltonian was developed in 1952 due to Dirac \cite{[2]} and is based on the  natural idea that if a sufficient number of edges are present in the graph on $n$ vertices  (by keeping the minimum degree at a fairly high level) then a Hamilton cycle will exist. \\

\noindent\textbf{Theorem A \cite{[2]}.} Every graph with $\delta\ge\frac{1}{2}n$ is hamiltonian.\\

A direct link between the number of edges and Hamilton cycles was established in 1959 due to Erd\"{o}s  and Gallai \cite{[3]}.

\noindent\textbf{Theorem B \cite{[3]}.} Every graph with $q\ge\frac{1}{2}(n^2-3n+5)$ is hamiltonian.\\

Recently it was proved a little surprising and, in fact, a contrary statement ensuring the existence of a Hamilton cycle if the number of edges is less than  $\delta^2+\delta$.\\

\noindent\textbf{Theorem C \cite{[5]}.} Every graph with $q\le\delta^2+\delta-1$ is hamiltonian.\\

In this paper we present the exact analog of  Theorem C for dominating cycles.\\

\noindent\textbf{Theorem 1.} Let $G$ be a 2-connected graph. If  
$$
q\le\left\{ 
\begin{array}{lll}
8 & \mbox{if} & \mbox{ }\delta=2, \\ \frac{3(\delta-1)(\delta+2)-1}{2} & \mbox{if} & \mbox{ }%
\delta\ge3, 
\end{array}
\right. 
$$
then each longest cycle in $G$ is a dominating cycle.\\

To show that Theorem 1 is sharp, suppose first that $\delta =2$. The graph $K_1+2K_2$ shows that the connectivity condition $\kappa\ge2$ in Theorem 1 can not be relaxed by replacing it with $\kappa\ge1$.   The graph with vertex set $\{v_1,v_2,...,v_8\}$ and edge set 
$$
\{v_1v_2,v_2v_3,v_3v_4,v_4v_5,v_5v_6,v_6v_1,v_1v_7,v_7v_8,v_8v_4\},
$$
shows that the size bound $q\le8$ can not be relaxed by replacing it with $q\le9$. Finally, the graph $K_2+3K_1$ shows that the conclusion "each longest cycle  in $G$ is a dominating cycle" can not be strengthened by replacing it with "$G$ is hamiltonian". Now let $\delta\ge3$. The graph $K_1+2K_\delta$   shows that the connectivity condition $\kappa\ge2$ in Theorem 1 can not be relaxed by replacing it with $\kappa\ge1$.  Further, the graph $K_2+3K_{\delta-1}$  shows that the size bound $q\le(3(\delta-1)(\delta+2)-1)/2$ can not be relaxed by replacing it with $q\le3(\delta-1)(\delta+2)/2$. Finally, the graph $K_\delta+(\delta+1)K_1$  shows that the main conclusion "each longest cycle  in $G$ is a dominating cycle" can not be strengthened by replacing it with "$G$ is hamiltonian". So, Theorem 1 is best possible in all respects. 

The following theorems are useful.\\

\noindent\textbf{Theorem D \cite{[2]}.} Every 2-connected graph either has a Hamilton cycle or has a cycle of length at least $2\delta$. \\

\noindent\textbf{Theorem E \cite{[4]}.} Let $G$ be a graph, $C$ a longest cycle in $G$ and $P$ a longest path in $G\backslash C$ of length $\overline{p}$. Then $|C|\ge(\overline{p}+2)(\delta-\overline{p})$.\\

\noindent\textbf{Theorem F \cite{[6]}.} Let $G$ be a graph on $n$ vertices and $d(x)+d(y)\ge n$ for each nonadjacent vertices $x,y$. Then $G$ is hamiltonian.

\section{Notations and preliminaries}

The set of vertices of a graph $G$ is denoted by $V(G)$ and the set of edges by $E(G)$. For $S$ a subset of $V(G)$, we denote by $G\backslash S$ the maximum subgraph of $G$ with vertex set $V(G)\backslash S$. For a subgraph $H$ of $G$ we use $G\backslash H$ short for $G\backslash V(H)$. The neighborhood of a vertex $x\in V(G)$ will be denoted by $N(x)$. Set $d(x)=|N(x)|$. Furthermore, for a subgraph $H$ of $G$ and $x\in V(G)$, we define $N_H(x)=N(x)\cap V(H)$ and  $d_H(x)=|N_H(x)|$.

A simple cycle (or just a cycle) $C$ of length $t$ is a sequence $v_1v_2...v_tv_1$ of distinct vertices $v_1,...,v_t$ with $v_iv_{i+1}\in E(G)$ for each $i\in \{1,...,t\}$, where $v_{t+1}=v_1$. When $t=2$, the cycle $C=v_1v_2v_1$ on two vertices $v_1, v_2$ coincides with the edge $v_1v_2$, and when $t=1$, the cycle $C=v_1$ coincides with the vertex $v_1$. So, all vertices and edges in a graph can be considered as cycles of lengths 1 and 2, respectively. A graph $G$ is hamiltonian if $G$ contains a Hamilton cycle, i.e. a cycle of length $n$.  

Paths and cycles in a graph $G$ are considered as subgraphs of $G$. If $Q$ is a path or a cycle, then the length of $Q$, denoted by $|Q|$, is $|E(Q)|$. We write $Q$ with a given orientation by $\overrightarrow{Q}$. For $x,y\in V(Q)$, we denote by $x\overrightarrow{Q}y$ the subpath of $Q$ in the chosen direction from $x$ to $y$. For $x\in V(Q)$, we denote the $h$-th successor and the $h$-th predecessor of $x$ on $\overrightarrow{Q}$ by $x^{+h}$ and $x^{-h}$, respectively. We abbreviate $x^{+1}$ and $x^{-1}$ by $x^+$ and $x^-$, respectively. \\ 

\noindent\textbf{Special definitions}. Let $G$ be a graph, $C$ a longest cycle in $G$ and $P=x\overrightarrow{P}y$ a longest path in $G\backslash C$ of length $\overline{p}\ge0$. Let $\xi_1,\xi_2,...,\xi_s$ be the elements of $N_C(x)\cup N_C(y)$ occuring on $C$ in a consecutive order and let 
$$
I_i=\xi_i\overrightarrow{C}\xi_{i+1},       \     I_i^\ast=\xi_i^+\overrightarrow{C}\xi_{i+1}^-   \   \  (i=1,2,...,s),
$$
where $\xi_{s+1}=\xi_1$. 

$(\ast1)$ We call  $I_1,I_2,...,I_s$ elementary segments on $C$ induced by $N_C(x)\cup N_C(y)$.

$(\ast2)$ We call a path $L=z\overrightarrow{L}w$ an intermediate path between elementary segments $I_a$ and $I_b$ if
$$
z\in V(I_a^\ast),   \    w\in V(I_b^\ast),    \    V(L)\cap V(C\cup P)=\{z,w\}.
$$

$(\ast3)$ Denote by $M(I_{i_1},I_{i_2},...,I_{i_t})$ the set of all intermediate paths between elementary segments  $I_{i_1},I_{i_2},...,I_{i_t}$.\\

\noindent\textbf{Lemma 1.} Let $G$ be a graph, $C$ a longest cycle in $G$ and $P=x\overrightarrow{P}y$ a longest path in $G\backslash C$ of length $\overline{p}\ge1$. If  $|N_C(x)|\ge2$, $|N_C(y)|\ge2$ and  $N_C(x)\not=N_C(y)$ then
$$
|C|\ge\left\{ 
\begin{array}{lll}
3\delta+\max\{\sigma_1, \sigma_2\}-1\ge3\delta & \mbox{if} & \mbox{ }\overline{p}=1, \\ \max\{2\overline{p}+8, 4\delta-2\overline{p}\} & \mbox{if} & \mbox{ }%
\overline{p}\ge2, 
\end{array}
\right. 
$$
where $\sigma_1=|N_C(x)\backslash N_C(y)|$ and $\sigma_2=|N_C(y)\backslash N_C(x)|$.\\

\noindent\textbf{Lemma 2.} Let $G$ be a graph, $C$ a longest cycle in $G$ and $P=x\overrightarrow{P}y$ a longest path in $G\backslash C$ of length $\overline{p}\ge0$.  If $N_C(x)=N_C(y)$ and $|N_C(x)|\ge2$ then for each  elementary segments $I_a$ and $I_b$ induced by $N_C(x)\cup N_C(y)$,

(a1) if $L$ is an intermediate path between  $I_a$ and $I_b$ then 

$$
|I_a|+|I_b|\ge2\overline{p}+2|L|+4,
$$

(a2) if $M(I_a,I_b)\subseteq E(G)$ and $|M(I_a,I_b)|=i$\ \  $(i\in\{1,2,3\})$ then 
$$
|I_a|+|I_b|\ge2\overline{p}+i+5.
$$

\noindent\textbf{Lemma 3.} Let $G$ be a graph, $S$  a cut set in $G$ and $H$ a connected component of $G\backslash S$ of order $h$. Then     
$$
q_H\ge\frac{h(2\delta-h+1)}{2},
$$
where $q_H=|\{xy\in E(G) : \{x,y\}\cap V(H)\not=\emptyset\}|$.                    \\

\noindent\textbf{Lemma 4.} Let $G$ be a 2-connected graph. If  $\delta\ge(n-2)/3$ then either 
$$
q\ge\left\{ 
\begin{array}{lll}
9 & \mbox{if} & \mbox{ }\delta=2, \\ \frac{3(\delta-1)(\delta+2)}{2} & \mbox{if} & \mbox{ }%
\delta\ge3, 
\end{array}
\right. 
$$
or each longest cycle in $G$ is a dominating cycle.\\

\section{Proofs}

\noindent\textbf{Proof of Lemma 1}. Put
$$
A_1=N_C(x)\backslash N_C(y),    \    A_2=N_C(y)\backslash N_C(x), \   M=N_C(x)\cap N_C(y).
$$
By the hypothesis,  $N_C(x)\not=N_C(y)$, implying that 
$$
\max \{|A_1|,|A_2|\}\ge1.
$$ 
Let $\xi_1,\xi_2,...,\xi_s$ be the elements of $N_C(x)\cup N_C(y)$ occuring on $C$ in a consecutive order. Put $I_i=\xi_i\overrightarrow{C}\xi_{i+1}$ $(i=1,2,...,s)$, where $\xi_{s+1}=\xi_1$. Clearly, $s=|A_1|+|A_2|+|M|$. Since $C$ is extreme, $|I_i|\ge2$ $(i=1,2,...,s)$. Moreover, if $\{\xi_i,\xi_{i+1}\}\cap M\not=\emptyset$ for some $i\in\{1,2,...,s\}$ then $|I_i|\ge\overline{p}+2$. In addition, if either $\xi_i\in A_1$, $\xi_{i+1}\in A_2$ or $\xi_i\in A_2$, $\xi_{i+1}\in A_1$ then again $|I_i|\ge\overline{p}+2$. \\

\textbf{Case 1}. $\overline{p}=1$.

\textbf{Case 1.1}. $|A_i|\ge1$ $(i=1,2)$.

It follows that among  $I_1,I_2,...,I_s$ there are  $|M|+2$ segments of length at least $\overline{p}+2$. Observing also that each of the remaining $s-(|M|+2)$ segments has a length at least 2, we get 
$$
|C|\ge(\overline{p}+2)(|M|+2)+2(s-|M|-2)
$$
$$
=3(|M|+2)+2(|A_1|+|A_2|-2)
$$
$$
=2|A_1|+2|A_2|+3|M|+2.
$$

Since $|A_1|=d(x)-|M|-1$ and  $|A_2|=d(y)-|M|-1$, we have
$$
|C|\ge2d(x)+2d(y)-|M|-2\ge3\delta+d(x)-|M|-2.
$$
Recalling that $d(x)=|M|+|A_1|+1$, we get 
$$
|C|\ge3\delta+|A_1|-1=3\delta+\sigma_1-1.
$$
Analogously, $|C|\ge3\delta+\sigma_2-1$. So, 
$$
|C|\ge3\delta+\max \{\sigma_1,\sigma_2\}-1\ge3\delta.
$$

\textbf{Case 1.2}. Either $|A_1|\ge1, |A_2|=0$ or $|A_1|=0, |A_2|\ge1$.

Assume w.l.o.g. that $|A_1|\ge1$ and $|A_2|=0$, i.e. $|N_C(y)|=|M|\ge2$ and $s=|A_1|+|M|$ . Hence, among $I_1,I_2,...,I_s$ there are $|M|+1$ segments of length at least $\overline{p}+2=3$. Taking into account that each of the remaining $s-(|M|+1)$ segments has a length at least 2 and $|M|+1=d(y)$, we get
$$
|C|\ge 3(|M|+1)+2(s-|M|-1)=3d(y)+2(|A_1|-1)
$$
$$
\ge3\delta+|A_1|-1=3\delta+\max\{\sigma_1,\sigma_2\}-1\ge3\delta.
$$

\textbf{Case 2}. $\overline{p}\ge2$.

We first prove that $|C|\ge2\overline{p}+8$. Since $|N_C(x)|\ge2$ and $|N_C(y)|\ge2$,  there are at least two segments among $I_1,I_2,...,I_s$ of length at least $\overline{p}+2$. If $|M|=0$ then clearly $s\ge4$ and 
$$
|C|\ge2(\overline{p}+2)+2(s-2)\ge2\overline{p}+8.
$$
Otherwise, since $\max\{|A_1|,|A_2|\}\ge1$, there are at least three elementary segments of length at least $\overline{p}+2$, i.e. 
$$
|C|\ge3(\overline{p}+2)\ge2\overline{p}+8.
$$
So, in any case, $|C|\ge 2\overline{p}+8$. 

To prove that $|C|\ge 4\delta-2\overline{p}$, we distinguish two main cases.\\

\textbf{Case 2.1}. $|A_i|\ge1$ $(i=1,2)$.

It follows that among  $I_1,I_2,...,I_s$ there are $|M|+2$ segments of length at least $\overline{p}+2$. Further, since each of the remaining $s-(|M|+2)$ segments has a length at least 2, we get 
$$
|C|\ge (\overline{p}+2)(|M|+2)+2(s-|M|-2)
$$
$$
=(\overline{p}-2)|M|+(2\overline{p}+4|M|+4)+2(|A_1|+|A_2|-2)
$$
$$
\ge2|A_1|+2|A_2|+4|M|+2\overline{p}.
$$
Observing also that 
$$
|A_1|+|M|+\overline{p}\ge d(x),  \quad   |A_2|+|M|+\overline{p}\ge d(y),  
$$
we have 

$$
2|A_1|+2|A_2|+4|M|+2\overline{p}
$$
$$
\ge 2d(x)+2d(y)-2\overline{p}\ge4\delta-2\overline{p},
$$
implying that $|C|\ge4\delta-2\overline{p}$.

 \textbf{Case 2.2}. Either $|A_1|\ge1, |A_2|=0$ or $|A_1|=0, |A_2|\ge1$.

Assume w.l.o.g. that $|A_1|\ge1$ and $|A_2|=0$, i.e.  $|N_C(y)|=|M|\ge2$ and $s=|A_1|+|M|$. It follows that among  $I_1,I_2,...,I_s$ there are  $|M|+1$ segments of length at least $\overline{p}+2$. Observing also that $|M|+\overline{p}\ge d(y)\ge\delta$, i.e. $2\overline{p}+4|M|\ge 4\delta-2\overline{p}$, we get 
$$
|C|\ge(\overline{p}+2)(|M|+1)\ge(\overline{p}-2)(|M|-1)+2\overline{p}+4|M|
$$
$$
\ge 2\overline{p}+4|M|\ge4\delta-2\overline{p}.       \quad            \quad                    \rule{7pt}{6pt} 
$$

\noindent\textbf{Proof of Lemma 2}. Let $\xi_1,\xi_2,...,\xi_s$ be the elements of $N_C(x)$ occuring on $C$ in a consecutive order. Put $I_i=\xi_i\overrightarrow{C}\xi_{i+1}$ $(i=1,2,...,s)$, where $\xi_{s+1}=\xi_1.$  To prove $(a1)$, let $L=z\overrightarrow{L}w$ be an intermediate path between elementary segments $I_a$ and $I_b$ with $z\in V(I_a^\ast)$ and   $w\in V(I_b^\ast)$.
Put
$$
|\xi_a\overrightarrow{C}z|=d_1,  \   |z\overrightarrow{C}\xi_{a+1}|=d_2,      \      |\xi_b\overrightarrow{C}w|=d_3,       \     |w\overrightarrow{C}\xi_{b+1}|=d_4,
$$
$$
C^\prime=\xi_ax\overrightarrow{P}y\xi_b\overleftarrow{C}z\overrightarrow{L}w\overrightarrow{C}\xi_a.
$$
Clearly, 
$$
|C^\prime|=|C|-d_1-d_3+|L|+|P|+2.           
$$
Since $C$ is extreme, we have $|C|\ge|C^\prime|$, implying that $d_1+d_3\ge\overline{p}+|L|+2$.  By a symmetric argument, $d_2+d_4\ge\overline{p}+|L|+2$. Hence 
$$
|I_a|+|I_b|=\sum_{i=1}^4d_i\ge2\overline{p}+2|L|+4.
$$

To proof  $(a2)$, let  $M(I_a,I_b)\subseteq E(G)$ and $|M(I_a,I_b)|=i$\ \ $(i\in \{1,2,3\})$.\\

\textbf{Case 1}. $i=1$.

It follows that $M(I_a,I_b)$ consists of a single intermediate edge $L=zw$. By (a1), 
$$
|I_a|+|I_b|\ge2\overline{p}+2|L|+4=2\overline{p}+6.
$$

\textbf{Case 2}. $i=2$.

It follows that $M(I_a,I_b)$ consists of two edges $e_1,e_2$. Put $e_1=z_1w_1$ and $e_2=z_2w_2$, where $\{z_1,z_2\}\subseteq V(I_a^\ast)$ and $\{w_1,w_2\}\subseteq V(I_b^\ast)$.\\

\textbf{Case 2.1}. $z_1\not=z_2$ and $w_1\not=w_2$.

Assume w.l.o.g. that $z_1$ and $z_2$ occur in this order on $I_a$. \\

\textbf{Case 2.1.1}. $w_2$ and  $w_1$ occur in this order on $I_b$.

Put
$$
|\xi_a\overrightarrow{C}z_1|=d_1,  \   |z_1\overrightarrow{C}z_2|=d_2,      \    |z_2\overrightarrow{C}\xi_{a+1}|=d_3, 
$$
$$
|\xi_b\overrightarrow{C}w_2|=d_4,       \       |w_2\overrightarrow{C}w_1|=d_5,         \        |w_1\overrightarrow{C}\xi_{b+1}|=d_6, 
$$
$$
C^{\prime}=\xi_a\overrightarrow{C}z_1w_1\overleftarrow{C}w_2z_2\overrightarrow{C}\xi_b x\overrightarrow{P}y\xi_{b+1}\overrightarrow{C}\xi_a.
$$
Clearly, 
$$
|C^{\prime}|=|C|-d_2-d_4-d_6+|\{e_1\}|+|\{e_2\}|+|P|+2
$$
$$
=|C|-d_2-d_4-d_6+\overline{p}+4.
$$
Since $C$ is extreme, $|C|\ge |C^{\prime}|$, implying that $d_2+d_4+d_6\ge \overline{p}+4$. By a symmetric argument, $d_1+d_3+d_5\ge\overline{p}+4$. Hence
$$
|I_a|+|I_b|= \sum_{i=1}^6d_i\ge2\overline{p}+8.
$$

\textbf{Case 2.1.2}. $w_1$ and $w_2$ occur in this order on $I_b$

Putting
$$
C^{\prime}=\xi_a\overrightarrow{C}z_1w_1\overrightarrow{C}w_2z_2\overrightarrow{C}\xi_b x\overrightarrow{P}y\xi_{b+1}\overrightarrow{C}\xi_a,
$$
we can argue as in Case 2.1.1. \\

\textbf{Case 2.2}. Either $z_1=z_2$, $w_1\not=w_2$ or $z_1\not=z_2$, $w_1=w_2$.

Assume w.l.o.g. that $z_1\not=z_2$, $w_1=w_2$ and $z_1, z_2$ occur  in this order on $I_a$. Put
$$
|\xi_a\overrightarrow{C}z_1|=d_1,  \   |z_1\overrightarrow{C}z_2|=d_2,      \    |z_2\overrightarrow{C}\xi_{a+1}|=d_3, 
$$
$$
|\xi_b\overrightarrow{C}w_1|=d_4,       \           |w_1\overrightarrow{C}\xi_{b+1}|=d_5, 
$$
$$
C^{\prime}=\xi_a x\overrightarrow{P}y\xi_b\overleftarrow{C}z_1w_1\overrightarrow{C}\xi_a,
$$
$$
C^{\prime\prime}=\xi_a\overrightarrow{C}z_2w_1\overleftarrow{C}\xi_{a+1}x\overrightarrow{P}y\xi_{b+1}\overrightarrow{C}\xi_a.
$$

Clearly, 
$$
|C^{\prime}|=|C|-d_1-d_4+|\{e_1\}|+|P|+2=|C|-d_1-d_4+\overline{p}+3,
$$
$$
|C^{\prime\prime}|=|C|-d_3-d_5+|\{e_2\}|+|P|+2=|C|-d_3-d_5+\overline{p}+3.
$$
Since $C$ is extreme, $|C|\ge |C^{\prime}|$ and $|C|\ge |C^{\prime\prime}|$, implying that 
$$
d_1+d_4\ge \overline{p}+3,  \    d_3+d_5\ge \overline{p}+3. 
$$
Hence, 
$$
|I_a|+|I_b|= \sum_{i=1}^5d_i\ge d_1+d_3+d_4+d_5+1\ge2\overline{p}+7.
$$

\textbf{Case 3}. $i=3$.

It follows that $M(I_a,I_b)$ consists of three edges $e_1,e_2,e_3$. Let $e_i=z_iw_i$ $(i=1,2,3)$, where $\{z_1,z_2,z_3\}\subseteq V(I_a^\ast)$ and $\{w_1,w_2,w_3\}\subseteq V(I_b^\ast)$. If there are two independent edges among $e_1,e_2,e_3$ then we can argue as in Case 2.1. Otherwise, we can assume w.l.o.g. that $w_1=w_2=w_3$ and $z_1,z_2,z_3$ occur in this order on $I_a$. Put
$$
|\xi_a\overrightarrow{C}z_1|=d_1,  \   |z_1\overrightarrow{C}z_2|=d_2,      \    |z_2\overrightarrow{C}z_3|=d_3,      
$$
$$
|z_3\overrightarrow{C}\xi_{a+1}|=d_4,      \    |\xi_b\overrightarrow{C}w_1|=d_5,       \           |w_1\overrightarrow{C}\xi_{b+1}|=d_6, 
$$
$$
C^{\prime}=\xi_a x\overrightarrow{P}y\xi_b\overleftarrow{C}z_1w_1\overrightarrow{C}\xi_a,
$$
$$
C^{\prime\prime}=\xi_a\overrightarrow{C}z_3w_1\overleftarrow{C}\xi_{a+1}x\overrightarrow{P}y\xi_{b+1}\overrightarrow{C}\xi_a.
$$

Clearly, 
$$
|C^{\prime}|=|C|-d_1-d_5+|\{e_1\}|+\overline{p}+2,
$$
$$
|C^{\prime\prime}|=|C|-d_4-d_6+|\{e_3\}|+\overline{p}+2.
$$
Since $C$ is extreme, we have $|C|\ge |C^{\prime}|$ and $|C|\ge |C^{\prime\prime}|$, implying that 
$$
d_1+d_5\ge \overline{p}+3,   \    d_4+d_6\ge \overline{p}+3. 
$$
Hence, 
$$
|I_a|+|I_b|= \sum_{i=1}^6d_i\ge d_1+d_4+d_5+d_6+2\ge2\overline{p}+8.     \quad      \quad    \rule{7pt}{6pt} 
$$

\noindent\textbf{Proof of Lemma 3}. Put 
$$
V(H)=\{v_1,...,v_h\},        \    |N(v_i)\cap S|=\beta_i     \    \    (i=1,...,h).
$$

Observing that $h\ge d(v_i)-\beta_i+1\ge \delta-\beta_i+1$ for each $i\in \{1,2,...,h\}$, we have $\beta_i\ge\delta-h+1$ $(i=1,2,...,h)$. Therefore,
$$
q_H=q(H)+\sum_{i=1}^h\beta_i=\frac{1}{2}\sum_{i=1}^hd_H(v_i)+\sum_{i=1}^h\beta_i
$$
$$
=\frac{1}{2}\sum_{i=1}^h(d_H(v_i)+\beta_i)+\frac{1}{2}\sum_{i=1}^h\beta_i =\frac{1}{2}\sum_{i=1}^hd(v_i)+\frac{1}{2}\sum_{i=1}^h(\delta-h+1)
$$
$$
 \ge\frac{1}{2}h\delta+\frac{1}{2}h(\delta-h+1)=\frac{h(2\delta-h+1)}{2}.        \quad          \quad      \rule{7pt}{6pt}                   
$$

\noindent\textbf{Proof of Lemma 4}. Let $C$ be a longest cycle in $G$ and $P=x_1\overrightarrow{P}x_2$ a longest path in $G\backslash C$ of length $\overline{p}$. If $|V(P)|\le1$ then $C$ is a dominating cycle and we are done. Let $|V(P)|\ge2$, that is $\overline{p}\ge1$. By the hypothesis, 
$$
|C|+\overline{p}+1\le n\le3\delta+2. \eqno{(1)}
$$
Let $\xi_1,\xi_2,...,\xi_s$ be the elements of $N_C(x_1)\cup N_C(x_2)$ occuring on $C$ in a consecutive order. Put
$$
I_i=\xi_i\overrightarrow{C}\xi_{i+1},       \     I_i^\ast=\xi_i^+\overrightarrow{C}\xi_{i+1}^-  \  \ (i=1,2,...,s),
$$
where $\xi_{s+1}=\xi_1.$   Let $Q$ be a longest path in $G$ with $Q=\xi\overrightarrow{Q}\eta$ and $V(Q)\cap V(C)=\{\xi,\eta\}$. Since $C$ is extreme, we have $|\xi\overrightarrow{C}\eta|\ge|Q|$ and  $|\eta\overrightarrow{C}\xi|\ge|Q|$, implying that 
$$
|C|=|\xi\overrightarrow{C}\eta|+|\eta\overrightarrow{C}\xi|\ge2|Q|.                            \eqno{(2)}
$$

\textbf{Case 1}. $\delta=2$.

 Since $\kappa\ge2$ and $\overline{p}\ge1$, we have $|Q|\ge3$. By (2),
$$
|C|=|y\overrightarrow{C}z|+|z\overrightarrow{C}y|\ge2|Q|\ge6,
$$
implying that $q\ge|C|+|Q|\ge9$.\\

\textbf{Case 2}. $\delta=3$.

If $n\ge10$ then 
$$
q\ge \frac{n\delta}{2}\ge15=\frac{3(\delta-1)(\delta+2)}{2}.
$$ 
Let 
$$
n\le9.  \eqno{(3)}
$$

\textbf{Case 2.1}. $\overline{p}=1$.

By (1) and (3),
$$
|C|\le n-\overline{p}-1\le 7.            \eqno{(4)}
$$
Since $\overline{p}=1$ and $\delta=3$, we have $|N_C(x_i)|\ge2$ $(i=1,2)$. If $N_C(x_1)\not=N_C(x_2)$ then by Lemma 1, $|C|\ge 3\delta=9$, contradicting (4). Let $N_C(x_1)=N_C(x_2)$. Further, since $C$ is extreme and $\overline{p}=1$, we have $|I_i|\ge3$ $(i=1,2,...,s)$. If $s\ge3$ then $|C|=\sum_{i=1}^s|I_i|\ge3s\ge9$, contradicting (4). Let $s=2$. If $M(I_1,I_2)\not=\emptyset$ then by Lemma 2, $|C|=|I_1|+|I_2|\ge 2\overline{p}+6=8$, contradicting (4). Thus, $M(I_1,I_2)=\emptyset$, implying that $G\backslash \{\xi_1,\xi_2\}$ is disconnected. Let $H_1,H_2,...,H_t$ be the connected components of $G\backslash \{\xi_1,\xi_2\}$. Clearly, $t\ge3$.  Put
$$
h_i=|V(H_i)|, \  q_i=|\{xy\in E(G) : \{x,y\}\cap V(H_i)\not=\emptyset\}|  \   (i=1,2,...,t).        \eqno{(5)}
$$
Assume w.l.o.g. that  $V(I_i^\ast)\subseteq V(H_i)$ $(i=1,2)$ and $V(H_3)=\{x_1,x_2\}$. It means that $h_i\ge2$  $(i=1,2,3)$. If $h_i\ge4$ for some $i\in\{1,2\}$ then 
$$
|C|\ge h_1+h_2+|\{\xi_1,\xi_2\}|\ge 8,
$$
contradicting (4). Let $2\le h_i\le3$ $(i=1,2,3)$. By Lemma 3, 
$$
q_i\ge\frac{h_i(2\delta-h_i+1)}{2}=\frac{h_i(7-h_i)}{2}\ge5   \quad   (i=1,2,3).       
$$
Hence 
$$
q\ge\sum_{i=1}^3q_i\ge15=\frac{3(\delta-1)(\delta+2)}{2}.
$$

\textbf{Case 2.2}. $\overline{p}\ge2$.

By (1) and (3), $|C|\le n-\overline{p}-1\le 6$.\\

\textbf{Case 2.2.1}. There is a cycle in $G\backslash C$.

Let $C^\prime$ be a cycle in $G\backslash C$. Since $\kappa\ge2$, there are two disjoint paths connecting $C^\prime$ and $C$, implying that $|Q|\ge4$. By (2), $|C|\ge 2|Q|\ge8$, contradicting (4).\\

\textbf{Case 2.2.2}. $G\backslash C$ is acyclic.

It follows that 
$$
|N_C(x_i)|\ge |N(x_i)|-1\ge\delta-1=2   \   \    (i=1,2).
$$
Hence $|Q|\ge\overline{p}+2\ge4$. By (2), $|C|\ge2|Q|\ge8$, contradicting (4). \\

\textbf{Case 3}. $\delta=4$.

If $n\ge14$ then 
$$
q\ge \frac{n\delta}{2}\ge28>\frac{3(\delta-1)(\delta+2)}{2}.
$$
Let 
$$
n\le13.  \eqno{(6)}
$$

\textbf{Case 3.1}. $\overline{p}=1$.

By (1) and (6),
$$
|C|\le n-\overline{p}-1\le11.                    \eqno{(7)}
$$
Since $\overline{p}=1$ and $\delta=4$, we have $|N_C(x_i)|\ge3$ $(i=1,2)$. If $N_C(x_1)\not=N_C(x_2)$ then by Lemma 1, $|C|\ge3\delta=12$, contradicting (7). Let $N_C(x_1)=N_C(x_2)$. Further, since $C$ is extreme and $\overline{p}=1$,  we have $|I_i|\ge3$ $(i=1,...,s)$. If $s\ge4$ then $|C|\ge3s\ge12$, contradicting (7). Thus $s=3$. \\

\textbf{Case 3.1.1}. $M(I_1,I_2,I_3)=\emptyset$.

It follows that $G\backslash \{\xi_1,\xi_2,\xi_3\}$ is disconnected. Let $H_1,H_2,...,H_t$ be the connected components of $G\backslash \{\xi_1,\xi_2,\xi_3\}$. Clearly, $t\ge4$. Assume w.l.o.g. that $V(I_i^\ast)\subseteq V(H_i)$ $(i=1,2,3)$ and $V(H_4)=\{x_1,x_2\}$. Using notation (5), we have $h_i\ge2$ $(i=1,2,3)$ and $h_4=2$. If $h_i\ge5$ for some $i\in \{1,2,3\}$ then clearly 
$$
|C|\ge\sum_{i=1}^3h_i+|\{\xi_1,\xi_2,\xi_3\}|\ge12,
$$
contradicting (7). Let $2\le h_i\le4$ $(i=1,2,3)$. By Lemma 3,
$$
q_i\ge \frac{h_i(2\delta-h_i+1)}{2}=\frac{h_i(9-h_i)}{2}\ge7  \quad (i=1,2,3,4).         
$$
So, 
$$
q\ge\sum_{i=1}^4q_i\ge 28>\frac{3(\delta-1)(\delta+2)}{2}.
$$

\textbf{Case 3.1.2}. $M(I_1,I_2,I_3)\not=\emptyset$.

Assume w.l.o.g. that $M(I_1,I_2)\not=\emptyset$, i.e. there is an intermediate path $L$ between $I_1$ and $I_2$.  By Lemma 2, 
$$
|I_1|+|I_2|\ge2\overline{p}+2|L|+4=2|L|+6.
$$
If $|L|\ge2$ then $|I_1|+|I_2|\ge10$ and hence $|C|=|I_1|+|I_2|+|I_3|\ge13$, contradicting (7). Otherwise,  $|L|=1$, implying that $M(I_1,I_2,I_3)\subseteq E(G)$. If $|M(I_1,I_2)|\ge2$ then by Lemma 2, $|I_1|+|I_2|\ge2\overline{p}+7=9$ and $|C|=\sum_{i=1}^3|I_i|\ge12$, contradicting (7). So, $|M(I_1,I_2)|=1$. By Lemma 2, $|I_1|+|I_2|\ge2\overline{p}+6=8$. Since $|I_3|\ge3$, we have $|C|=\sum_{i=1}^3|I_i|\ge11$. By (1), $n\ge|C|+\overline{p}+1\ge13$. Combining  $n\ge13$ and $|C|\ge11$ with (6) and (7), we get
$$
n=13, \ |C|=11, \ |I_1|+|I_2|=8, \  |I_3|=3, \  V(G)=V(C\cup P).   \eqno{(8)}
$$
Since $|I_1|+|I_2|=8$ and $|I_i|\ge3$ $(i=1,2)$, we can assume w.l.o.g. that either $|I_1|=|I_2|=4$ or $|I_1|=3$, $|I_2|=5$. If $|I_1|=|I_2|=4$ then by Lemma 2, $M(I_1,I_3)=M(I_2,I_3)=\emptyset$, implying that $|M(I_1,I_2,I_3)|=1$. Further, if $|I_1|=3$ and $|I_2|=5$ then by Lemma 2, $M(I_1,I_3)=\emptyset$ and $|M(I_2,I_3)|\le1$, implying that $|M(I_1,I_2,I_3)|\le2$. So, in any case, 
$$
1\le |M(I_1,I_2,I_3)|\le2.         \eqno{(9)}
$$
Let $e\in M(I_1,I_2,I_3)$ and  $e=zw$. Put $G^\prime=G\backslash e$. Form a graph $G^{\prime\prime}$ in the following way. If $d(z)\ge\delta=4$ and $d(w)\ge\delta=4$ in $G^\prime$ then we take $G^{\prime\prime}=G^\prime$. Next, we let $d(z)=\delta-1=3$ and $d(w)\ge\delta=4$ in $G^\prime$. If $\{\xi_1,\xi_2,\xi_3\}\subseteq N(z)$ then clearly $d(z)\ge4$ in $G^\prime$, contradicting the hypothesis. Otherwise, $zv\not\in E(G^\prime)$ for some $v\in \{\xi_1,\xi_2,\xi_3\}$ and we take $G^{\prime\prime}=G^\prime+ \{zv\}$. Finally, if $d(z)=d(w)=3$ then as above, $zv\not\in E(G^\prime)$ and $wu\not\in E(G^\prime)$ for some $v,u\in \{\xi_1,\xi_2,\xi_3\}$ and we take $G^{\prime\prime}=G^\prime+ \{zv,wu\}$. Clearly, $\delta(G^{\prime\prime})=\delta(G)=4$ and $q=q(G)\ge q(G^{\prime\prime})-1$. Furthermore, deleting step by step all edges from $M(I_1,I_2,I_3)$ and adding at most two appropriate new edges against each deleting edge, we can form a graph $G^\ast$ with $\delta(G^\ast)=\delta(G)=4$ and $q(G)\ge q(G^\ast)-|M(I_1,I_2,I_3)|$. By (9), $q(G)\ge q(G^\ast)-2$. In fact, $G^\ast= (G\backslash M(I_1,I_2,I_3))+ E^\ast$, where  $E^\ast$ consists of at most $2|M(I_1,I_2,I_3)|$ appropriate new edges having exactly one end in common with $\{\xi_1,\xi_2,\xi_3\}$, implying  that $G^\ast\backslash \{\xi_1,\xi_2,\xi_3\}$ is disconnected. Let $H_1,H_2,H_3,H_4$ be the connected components of $G^\ast\backslash \{\xi_1,\xi_2,\xi_3\}$ with $V(H_i)=V(I_i^\ast)$ $(i=1,2,3)$ and $V(H_4)=\{x_1,x_2\}$. Put
$$
h_i=|V(H_i)|, \  q_i=|\{xy\in E(G^\ast) : \{x,y\}\cap V(H_i)\not=\emptyset\}|  \  \ (i=1,2,3,4).       
$$
Since $|I_1|+|I_2|=8$, we have either $|I_1|\ge4$ or $|I_2|\ge4$. Assume w.l.o.g. that $|I_1|\ge4$, that is $h_1\ge3$. As in Case 3.1.1, $2\le h_i\le4$ $(i=1,2,3,4)$. By Lemma 3, 
$$
q_1(G^\ast)\ge\frac{h_1(2\delta-h_1+1)}{2}=\frac{h_1(9-h_1)}{2}\ge9,
$$
$$
q_i(G^\ast)\ge\frac{h_i(2\delta-h_i+1)}{2}=\frac{h_i(9-h_i)}{2}\ge7 \quad  (i=2,3,4).
$$
Hence $q(G^\ast)\ge \sum_{i=1}^4q_i(G^\ast)\ge30$, implying that 
$$
q(G)\ge q(G^\ast)-2\ge28>\frac{3(\delta-1)(\delta+2)}{2}.
$$

\textbf{Case 3.2}. $\overline{p}=2$.

Put $P=x_1x_3x_2$. By (1) and (6),
$$
|C|\le n-\overline{p}-1\le10.          \eqno{(10)}
$$
Since $\delta=4$ and $\overline{p}=2$, we have $|N_C(x_i)|\ge2$ $(i=1,2)$. If $N_C(x_1)\not=N_C(x_2)$ then by Lemma 1, $|C|\ge4\delta-2\overline{p}=12$, contradicting (10). Let $N_C(x_1)=N_C(x_2)$. Recalling that $C$ is extreme and $\overline{p}=2$, we conclude that $|I_i|\ge4$ $(i=1,2,...,s)$. If $s\ge3$ then $|C|\ge4s\ge12$, contradicting (10). Let $s=2$, implying that $x_1x_2\in E(G)$. By symmetric arguments, we can state that $N_C(x_1)=N_C(x_2)=N_C(x_3)$. \\

\textbf{Case 3.2.1}. $M(I_1,I_2)=\emptyset$.

Let $H_1,H_2,...,H_t$ be the connected components of $G\backslash \{\xi_1,\xi_2\}$.  Assume w.l.o.g. that $V(I_i^\ast)\subseteq V(H_i)$ $(i=1,2)$ and $V(H_3)=\{x_1,x_2,x_3\}$. Using notation (5), we have $h_i\ge3$ $(i=1,2,3)$. If $h_i\ge6$ for some $i\in \{1,2\}$, then  $|C|\ge  h_1+h_2+|\{\xi_1,\xi_2\}|\ge11$, contradicting (10). Let $3\le h_i\le5$ $(i=1,2,3)$. By Lemma 3, 
$$
q_i\ge\frac{h_i(2\delta-h_i+1)}{2}=\frac{h_i(9-h_i)}{2}\ge9 \quad  (i=1,2,3).
$$ 
Hence 
$$
q\ge \sum_{i=1}^3q_i\ge 27=\frac{3(\delta-1)(\delta+2)}{2}.
$$

\textbf{Case 3.2.2}. $M(I_1,I_2)\not=\emptyset$.

By the definition, there is an intermediate path $L$ between $I_1$ and $I_2$.  By Lemma 2, 
$$
|I_1|+|I_2|\ge2\overline{p}+2|L|+4=2|L|+8.
$$
If $|L|\ge2$ then $|C|=|I_1|+|I_2|\ge12$, contradicting (10). Otherwise,  $|L|=1$ and therefore,  $M(I_1,I_2)\subseteq E(G)$. If $|M(I_1,I_2)|\ge2$ then by Lemma 2, 
$$
|C|=|I_1|+|I_2|\ge 2\overline{p}+7=11,
$$
contradicting (10). Now let $|M(I_1,I_2)|=1$, i.e. $M(I_1,I_2)$ consists of a single edge $e$. By Lemma 2, 
$$
|C|=|I_1|+|I_2|\ge 2\overline{p}+6=10,
$$
and by (1), $n\ge |C|+\overline{p}+1\ge13$. Combining $n\ge13$ and $|C|\ge10$ with (6) and (10), we get
$$
|C|=|I_1|+|I_2|=10,       \     n=13,    \      V(G)=V(C\cup P).        \eqno{(11)}
$$
Put $G^\prime =G\backslash e$ and let $H_1,H_2,H_3$ be the connected components of $G^\prime\backslash \{\xi_1,\xi_2\}$ with $V(H_i)=V(I_i^\ast)$ $(i=1,2)$ and $V(H_3)=V(P)$. Since $|I_1|+|I_2|=10$, we can assume w.l.o.g. that $|I_1|\ge5$. Using notation (5) for $G^\prime$, we have $h_1\ge|I_1|-1\ge4$ and $h_i\ge3$ $(i=2,3)$. If $h_i\ge6$ for some $i\in \{1,2\}$ then $|C|\ge h_1+h_2+|\{\xi_1,\xi_2\}|\ge11$, contradicting (11). Let $4\le h_1\le5$ and $3\le h_i\le5$ $(i=2,3)$. If $\delta(G^\prime)=\delta(G)$ then we can argue as in Case 3.2.1. Otherwise, as in Case 3.1.2, we can form a  graph $G^\ast$ by adding at most two new edges in $G^\prime$ such that $\delta(G^\ast)=\delta(G)$ and $G^\ast\backslash \{\xi_1,\xi_2\}$ has exactly three connected components. Recalling that $4\le h_1\le 5$ and using Lemma 3, we get
$$
q_1(G^\ast)\ge\frac{h_1(2\delta-h_1+1)}{2}=\frac{h_1(9-h_1)}{2}=10,
$$
$$
q_i(G^\ast)\ge\frac{h_i(2\delta-h_i+1)}{2}=\frac{h_i(9-h_i)}{2}\ge9 \quad  (i=2,3).
$$
So, $q(G^\ast)\ge\sum_{i=1}^3q_i(G^\ast)\ge28$, implying that 
$$
q(G)\geq(G^\ast)-1\ge27=\frac{3(\delta-1)(\delta+2)}{2}.
$$

\textbf{Case 3.3}. $\overline{p}=3$.

By (1) and (6),
$$
|C|\le n-\overline{p}-1\le9.          \eqno{(12)}
$$
Since $\delta=4$ and $\overline{p}=3$, we have $|N_C(x_i)|\ge1$ $(i=1,2)$. If  $|N_C(x_i)|\ge2$ for some $i\in \{1,2\}$ then $|Q|\ge \overline{p}+2=5$ and by (2), $|C|\ge2|Q|\ge10$, contradicting (12). Let $|N_C(x_i)|=1$ $(i=1,2)$. If $N_C(x_1)\not=N_C(x_2)$ then again $|Q|\ge5$ and $|C|\ge10$, contradicting (12). Thus, $N_C(x_1)=N_C(x_2)$. It follows that $G[V(P)]$ is complete. Since $\kappa\ge2$, there are two disjoint paths connecting $G[V(P)]$ and $C$, implying that $|Q|\ge5$ and $|C|\ge10$, contradicting (12).\\

\textbf{Case 3.4}. $\overline{p}=4$.

Put $P=x_1x_3x_4x_5x_2$. By (1) and (6),
$$
|C|\le n-\overline{p}-1\le8.          \eqno{(13)}
$$

\textbf{Case 3.4.1}. $x_1x_2\in E(G)$.

Put $C^\prime =x_1x_2x_5x_4x_3x_1$. Since $\kappa\ge2$, there are two disjoint paths connecting $C^\prime$ and $C$. Since $|C^\prime|=5$, we have $|Q|\ge5$ and by (2), $|C|\ge2|Q|\ge10$, contradicting (13). \\

 \textbf{Case 3.4.2}. $x_1x_2\not\in E(G)$.

As in Case 3.3, it can be shown that 
$$
N_C(x_1)=N_C(x_2),          \  \       |N_C(x_1)|=|N_C(x_2)|=1.
$$
Since $\delta=4$, we have $\{x_1x_4,x_1x_5,x_2x_3,x_2x_4\}\subset E(G)$. Hence, $x_1x_4x_5x_2x_3x_1$ is a Hamilton cycle  in $G[V(P)]$ and we can argue as in Case 3.4.1.\\

\textbf{Case 3.5}. $\overline{p}\ge5$.

If $|C|=n$ then $C$ is a dominating cycle. Otherwise, by Theorem D, $|C|\ge2\delta=8$. On the other hand, by (1) and (6), $|C|\le n-\overline{p}-1\le7$, a contradiction.\\

\textbf{Case 4}. $\delta\ge5$.

If $C$ is a Hamilton cycle then we are done. Otherwise, by Theorem D,   
$$
c\ge2\delta.                           \eqno{(14)}
$$
By (1), $\overline{p}\le 3\delta-|C|+1\le\delta+1$. So, 
$$
1\le\overline{p}\le\delta+1.     \eqno{(15)}
$$
We  distinguish two main cases, namely $1\le\overline{p}\le \delta-3$ and $\delta-2\le\overline{p}\le\delta+1$.\\

\textbf{Case 4.1}. $1\le\overline{p}\le \delta-3$.

It follows that 
$$
(\overline{p}+2)(\delta-\overline{p})=(\overline{p}-1)(\delta-\overline{p}-3)+3\delta-3\ge3\delta-3.
$$
By Theorem E,
$$
|C|\ge3\delta-3.                        \eqno{(16)}
$$
By (1) and (16), $\overline{p}\le 3\delta-|C|+1\le4$. So,
$$
1\le\overline{p}\le4.               \eqno{(17)}
$$

\textbf{Case 4.1.1}. $\overline{p}=1$.

It follows that  $|N_C(x_i)|\ge\delta-1>2$ $(i=1,2)$. By (1), $|C|\le 3\delta+1-\overline{p}=3\delta$. Combining this with (16), we have
$$
3\delta-3\le |C|\le3\delta.                  
$$

\textbf{Case 4.1.1.1}. $3\delta-3\le|C|\le3\delta-1$. 

If $N_C(x_1)\not=N_C(x_2)$ then by Lemma 1, $|C|\ge 3\delta$, contradicting the hypothesis. Let $N_C(x_1)=N_C(x_2)$. Since $C$ is extreme and $\overline{p}=1$,  we have $|I_i|\ge3$ $(i=1,...,s)$. If $s\ge \delta$ then $|C|\ge3s\ge3\delta$, again contradicting the hypothesis. Let $s\le\delta-1$. On the other hand, $s=|N_C(x_1)|=d(x_1)-1\ge\delta-1$, implying that $s=\delta-1$. \\

\textbf{Claim 1}. $M(I_1,I_2,...,I_s)\subseteq E(G)$ and $|M(I_1,I_2,...,I_s)|\le\delta-2$.

Assume first that $3\delta-3\le|C|\le3\delta-2$. If $M(I_a,I_b)\not=\emptyset$ for some two elementary segments $I_a$ and $I_b$ then by Lemma 2, $|I_a|+|I_b|\ge2\overline{p}+6=8$, implying that $|C|\ge 3\delta-1$, a contradiction. Otherwise, $|M(I_1,I_2,...,I_s)|=0<\delta-2$. Now let $|C|=3\delta-1$. If $M(I_1,I_2,...,I_s)=\emptyset$ then we are done. Let $M(I_1,I_2,...,I_s)\not=\emptyset$, i.e.   $M(I_a,I_b)\not=\emptyset$ for some elementary segments $I_a$ and $I_b$. By the definition, there is an intermediate path $L$ between $I_a$ and $I_b$. If $|L|\ge2$ then by Lemma 2, $|I_a|+|I_b|\ge2\overline{p}+2|L|+4=10$, implying that $|C|\ge 3\delta$, a contradiction. Otherwise, $M(I_1,I_2,...,I_s)\subseteq E(G)$ and by Lemma 2, $|I_a|+|I_b|\ge2\overline{p}+6=8$, i.e.  $|C|\ge3\delta-1$. Recalling that $|C|=3\delta-1$, we can state that 
$$
|I_a|+|I_b|=8 \  \mbox{and} \   |I_i|=3   \   \mbox{for each}   \    i\in \{1,2,...,s\}\backslash \{a,b\}.
$$ 
If $|I_a|=|I_b|=4$ then by Lemma 2, $M(I_i,I_j)=\emptyset$ if $\{i,j\}\not=\{a,b\}$, i.e. $|M(I_1,I_2,...,I_s)|=1<\delta-2$. Otherwise, assume w.l.o.g. that $|I_a|=5$ and $|I_b|=3$, i.e. $|I_a|=5$ and $|I_i|=3$ for each $i\in \{1,2,...,s\}\backslash \{a\}$. As above,  $|M(I_a,I_i)|\le1$ for each $i\in \{1,2,...,s\}\backslash \{a\}$. Observing also that $M(I_i,I_j)=\emptyset$ for each distinct $i,j$ if $a\not\in \{i,j\}$, we conclude that $|M(I_1,I_2,...,I_s)|\le s-1= \delta-2$. Claim 1 is proved.    \quad     $\Delta$         \\

Put $G^\prime=G\backslash M(I_1,I_2,...,I_s)$. As in Case 3.1.2, we can form a graph $G^\ast$ by adding at most $2|M(I_1,I_2,...,I_s)|$ new edges in $G^\prime$ such that $\delta(G^\ast)=\delta(G)$, $G^\ast\backslash \{\xi_1,\xi_2,...,\xi_s\}$ is disconnected and $q(G)\ge q(G^\ast)-|M(I_1,I_2,...,I_s)|$. By Claim 1, 
$$
q(G)\ge q(G^\ast)-\delta+2.                \eqno{(18)}
$$
Let $H_1,H_2,...,H_{s+1}$ be the connected components of $G^\ast\backslash \{\xi_1,\xi_2,...,\xi_s\}$ with $V(I_i^\ast)\subseteq V(H_i)$ $(i=1,2,...,s)$ and $V(H_{s+1})=\{x_1,x_2\}$. Using notation (5) for $G^\ast$, we have $h_i\ge2$ $(i=1,2,...,s+1)$. If $h_i\ge6$ for some $i\in \{1,2,...,s\}$ then $n\ge3\delta+3$, contradicting (1). Let $2\le h_i\le 5<2\delta-1$ $(i=1,2,...,s+1)$. It follows that $(h_i-2)(2\delta-h_i-1)\ge0$ which is equivalent to
$$
\frac{h_i(2\delta-h_i+1)}{2}\ge2\delta-1 \quad    (i=1,2,...,s+1).
$$
By Lemma 3, $q_i(G^\ast)\ge2\delta-1$ $(i=1,2,...,s+1)$, implying that
$$
q(G^\ast)\ge\sum_{i=1}^{s+1}q_i(G^\ast)\ge(s+1)(2\delta-1)=\delta(2\delta-1).
$$
By (18),
$$
q\geq(G^\ast)-\delta+2\ge2(\delta^2-\delta+1)\ge\frac{3(\delta-1)(\delta+2)}{2}.
$$

\textbf{Case 4.1.1.2}. $|C|=3\delta$.

\textbf{Case 4.1.1.2.1}. $N_C(x_1)\not= N_C(x_2)$.

It follows that $\max \{\sigma_1,\sigma_2\}\ge1$, where 
$$
\sigma_1=|N_C(x_1)\backslash N_C(x_2)|,   \quad \sigma_2=|N_C(x_2)\backslash N_C(x_1)|.
$$
If $\max \{\sigma_1,\sigma_2\}\ge2$ then by Lemma 1, $|C|\ge3\delta+1$, contradicting the hypothesis. Let $\max \{\sigma_1,\sigma_2\}=1$. Clearly $s\ge\delta$ and $|I_i|\ge3$ $(i=1,2,...,s)$. If $s\ge \delta+1$ then $|C|\ge3s\ge3\delta+3$, a contradiction. Let $s=\delta$, implying that $|I_i|=3$ $(i=1,2,...,s)$. By Lemma 2, $M(I_1,I_2,...,I_s)=\emptyset$. Let $H_1,H_2,...,H_{s+1}$ be the connected components of $G\backslash \{\xi_1,\xi_2,...,\xi_s\}$ with $V(H_i)=V(I^\ast_i)$ $(i=1,2,...,s)$ and $V(H_{s+1})=\{x_1,x_2\}$. Using notation (5), we have $h_i=2$ $(i=1,2,...,s+1)$. By Lemma 3,
$$
q_i\ge\frac{h_i(2\delta-h_i+1)}{2}=2\delta-1  \quad   (i=1,2,...,s+1),
$$
implying that
$$
q\ge \sum^{s+1}_{i=1}q_i\ge (s+1)(2\delta-1)=(\delta+1)(2\delta-1)>\frac{3(\delta-1)(\delta+2)}{2}.
$$

\textbf{Case 4.1.1.2.2}. $N_C(x_1)= N_C(x_2)$.

Clearly, $s\ge\delta-1$. If $s\ge\delta$ then we can argue as in Case 4.1.1.2.1. Let $s=\delta-1$. If $|I_i|+|I_j|\ge10$ for some distinct $i,j\in \{1,2,...,s\}$ then $|C|\ge 10+3(s-2)=3\delta+1$, contradicting the hypothesis. Hence 
$$
|I_i|+|I_j|\le9 \quad  \mbox{for each distinct} \quad   i,j\in \{1,2,...,s\}.           \eqno{(19)}
$$

\textbf{Claim 2}. $M(I_1,I_2,...,I_s)\subseteq E(G)$ and

$(\ast1)$ if $\max_i|I_i|\le4$ then $|M(I_1,I_2,...,I_s)|\le3$,

$(\ast2)$ if $\max_i|I_i|=5$ then $|M(I_1,I_2,...,I_s)|\le\delta-1$,

$(\ast3)$ if $\max_i|I_i|=6$ then $|M(I_1,I_2,...,I_s)|\le2(\delta-2)$.

\textbf{Proof}. If $M(I_1,I_2,...,I_s)=\emptyset$ then we are done. Otherwise, $M(I_a,I_b)\not=\emptyset$ for some distinct $a,b\in \{1,2,...,s\}$. By the definition, there is an intermediate path $L$ between $I_a$ and $I_b$. If $|L|\ge2$ then by Lemma 2, 
$$
|I_a|+|I_b|\ge2\overline{p}+2|L|+4\ge10,
$$
contradicting (19). Otherwise, $|L|=1$ and $M(I_1,I_2,...,I_s)\subseteq E(G)$. By Lemma 2, $|I_a|+|I_b|\ge2\overline{p}+6=8$. Combining this with (19), we have
$$
8\le |I_a|+|I_b|\le9.               \eqno{(20)}
$$
Furthermore, if $|M(I_a,I_b)|\ge3$ then by Lemma 2, $|I_a|+|I_b|\ge2\overline{p}+8=10$,  contradicting (20). So,
$$
1\le |M(I_i,I_j)|\le2   \   \mbox{for each distinct}  \   i,j\in \{1,2,...,s\}.
$$
Put $r=|\{i\mid |I_i|\ge4\}|$. If $r\ge4$ then $|C|\ge3(s-4)+16=3\delta+1$, contradicting the hypothesis. Further, if $r=0$ then by Lemma 2, $M(I_1,I_2,...,I_s)=\emptyset$. Let $1\le r\le3$. \\

\textbf{Case a1}. $r=3$.

It follows that $|I_{a_i}|\ge4$ $(i=1,2,3)$ for some distinct $a_1,a_2,a_3\in\{1,2,...,s\}$ and $|I_i|=3$ for each $i\in \{1,2,...,s\}\backslash \{a_1,a_2,a_3\}$. Since $s=\delta-1$ and $|C|=3\delta$, we have $|I_{a_1}|=|I_{a_2}|=|I_{a_3}|=4$, i.e. $\max|I_i|=4$. By Lemma 2, $|M(I_{a_i},I_{a_j})|\le1$ for each distinct $i,j\in\{1,2,3\}$. Moreover, we have $|M(I_i,I_j)|=0$ if either $i\not\in \{i_1,i_2,i_3\}$ or $j\not\in \{i_1,i_2,i_3\}$. So, $|M(I_1,I_2,...,I_s)|\le3$.\\

\textbf{Case a2}. $r=2$.

It follows that $|I_a|\ge4$ and $|I_b|\ge4$ for some $a,b\in \{1,2,...,s\}$ and $|I_i|=3$ for each $i\in \{1,2,...,s\}\backslash \{a,b\}$. By (20), we can assume w.l.o.g. that either $|I_a|=|I_b|=4$ or $|I_a|=5$, $|I_b|=4$.\\

\textbf{Case a2.1}. $|I_a|=|I_b|=4$.

It follows that $\max_i |I_i|=4$. By Lemma 2, $|M(I_a,I_b)|\le1$ and $M(I_i,I_j)=\emptyset$ if $\{i,j\}\not=\{a,b\}$, implying that $|M(I_1,I_2,...,I_s)|\le1$.\\

\textbf{Case a2.2}. $|I_a|=5$, $|I_b|=4$.

It follows that $\max_i |I_i|=5$. By Lemma 2, we have $|M(I_a,I_b)|\le2$ and $|M(I_a,I_i)|\le1$ for each $i\in \{1,2,...,s\}\backslash \{a,b\}$ and $M(I_i,I_j)=\emptyset$ if $a\not\in\{i,j\}$. 
Thus, $|M(I_1,I_2,...,I_s)|\le\delta-1$ .\\

\textbf{Case a3}. $r=1$.

It follows that $|I_a|\ge4$ for some $a\in \{1,2,...,s\}$ and $|I_i|=3$ for each $i\in \{1,2,...,s\}\backslash \{a\}$. By (20), $4\le|I_a|\le6$.\\

\textbf{Case a3.1}. $|I_a|=4$.

It follows that $\max_1 |I_i|=4$. B Lemma 2, $M(I_a,I_i)=\emptyset$ for each  $i\in \{1,2,...,s\}\backslash \{a\}$, implying that $|M(I_1,I_2,...,I_s)|=0$. \\

\textbf{Case a3.2}. $|I_a|=5$.

It follows that $\max_i |I_i|=5$. B Lemma 2, $|M(I_a,I_i)|\le1$ for each  $i\in \{1,2,...,s\}\backslash \{a\}$ and $M(I_i,I_j)=\emptyset$ if $a\not\in\{i,j\}$, that is $|M(I_1,I_2,...,I_s)|\le\delta-2$. \\

\textbf{Case a3.3}. $|I_a|=6$.

It follows that $\max_i |I_i|=6$. By Lemma 2, $|M(I_a,I_i)|\le2$ for each  $i\in \{1,2,...,s\}\backslash \{a\}$ and $M(I_i,I_j)=\emptyset$ if $a\not\in\{i,j\}$, that is $|M(I_1,I_2,...,I_s)|\le2(\delta-2)$. Claim 2 is proved.   \quad        $\Delta$                \\

Put $G^\prime=G\backslash M(I_1,I_2,...,I_s)$. As in Case 3.1.2, we can form a graph $G^\ast$ by adding in $G^\prime$ at most $2|M(I_1,I_2,...,I_s)|$ new edges such that $\delta(G^\ast)=\delta(G)$, $G^\ast\backslash \{\xi_1,\xi_2,...,\xi_s\}$ is disconnected and 
$$
q(G)\ge q(G^\ast)-|M(I_1,I_2,...,I_s)|.   \eqno{(21)}
$$
Let $H_1,H_2,...,H_{s+1}$ be the connected components of $G^\ast\backslash \{\xi_1,\xi_2,...,\xi_s\}$ with $V(I_i^\ast)\subseteq V(H_i)$ $(i=1,2,...,s)$ and $V(H_{s+1})=\{x_1,x_2\}$. Using notation (5) for $G^\ast$, we have $h_i\ge2$ $(i=1,2,...,s+1)$. If $h_i\ge6$ for some $i\in \{1,2,...,s\}$ then $n\ge3\delta+3$, contradicting (1). Let $2\le h_i\le 5<2\delta-1$ $(i=1,2,...,s+1)$. It follows that 
$(h_i-2)(2\delta-h_i-1)\ge0$ which is equivalent to
$$
\frac{h_i(2\delta-h_i+1)}{2}\ge2\delta-1   \quad   (i=1,2,...,s+1).          \eqno{(22)}
$$

\textbf{Case 4.1.1.2.2.1}. $\max_i |I_i|\le4$.

By (22) and Lemma 3, $q_i(G^\ast)\ge2\delta-1$ $(i=1,2,...,s+1)$. Hence 
$$
q(G^\ast)\ge\sum_{i=1}^{s+1}q_i(G^\ast)\ge(s+1)(2\delta-1)=\delta(2\delta-1).
$$
Using (21) and Claim 2, we have
$$
q\ge q(G^\ast)-3\ge\delta(2\delta-1)-3\ge\frac{3(\delta-1)(\delta+2)}{2}.
$$

\textbf{Case 4.1.1.2.2.2}. $\max_i |I_i|=5$.

Assume w.l.o.g. that $\max_i |I_i|=|I_1|=5$, i.e. $4\le h_1\le5$. By (22) and Lemma 3, $q_i(G^\ast)\ge2\delta-1$ $(i=2,...,s+1)$ and 
$$
q_1(G^\ast)\ge\frac{h_1(2\delta-h_1+1)}{2}\ge2(2\delta-3).
$$
Hence 
$$
q(G^\ast)\ge s(2\delta-1)+2(2\delta-3)=2\delta^2+\delta-5.
$$
By (21) and Claim 2,
$$
q\ge q(G^\ast)-(\delta-1)\ge2\delta^2-4>\frac{3(\delta-1)(\delta+2)}{2}.
$$

\textbf{Case 4.1.1.2.2.3}. $\max_i|I_i|=6$.

Assume w.l.o.g. that $\max_i|I_i|=|I_1|=6$, that is $h_1=5$. By (22) and Lemma 3, $q_i(G^\ast)\ge 2\delta-1$ $(i=2,...,s+1)$ and 
$$
q_1(G^\ast)\ge\frac{h_1(2\delta-h_1+1)}{2}=5(\delta-2).
$$
Hence  
$$
q(G^\ast)\ge s(2\delta-1)+5(\delta-2)=2\delta^2+2\delta-9.
$$
By (21) and Claim 2,
$$
q\ge q(G^\ast)-2(\delta-2)\ge2\delta^2-5>\frac{3(\delta-1)(\delta+2)}{2}.
$$

\textbf{Case 4.1.2}. $\overline{p}=2$.

Put $P=x_1x_3x_2$. It follows that $|N_C(x_i)|\ge \delta-2>2$ $(i=1,2)$. By (1), $|C|\le3\delta+1-\overline{p}=3\delta-1$. Combining this with (16), we have
$$
3\delta-3\le |C|\le 3\delta-1.            \eqno{(23)}
$$
If $N_C(x_1)\not=N_C(x_2)$ then by Lemma 1, $|C|\ge4\delta-2\overline{p}=4\delta-4$. By (23), $4\delta-4\le |C|\le 3\delta-1$, a contradiction. Let $N_C(x_1)=N_C(x_2)$.  Since $C$ is extreme, we have $|I_i|\ge4$ $(i=1,2,...,s)$. If $s\ge \delta-1$ then $|C|\ge 4s\ge 4\delta-4\ge3\delta$, contradicting (23).  Hence $s\le \delta-2$. Recalling also that $s=|N_C(x_1)|\ge \delta-2$, we get $s=\delta-2$. It follows that $x_1x_2\in E(G)$. By a symmetric argument, $N_C(x_i)=N_C(x_1)$ $(i=2,3)$.\\

\textbf{Claim 3}. $|M(I_1,I_2,...,I_s)|\le 1$.

\textbf{Proof}. If $M(I_1,I_2,...,I_s)=\emptyset$ then we are done. Otherwise, $|M(I_a,I_b)|\ge1$ for some $a,b\in\{1,2,...,s\}$, i.e. there is an intermediate path $L$ between $I_a$ and $I_b$. If $|L|\ge2$ then by Lemma 2, 
$$
|I_a|+|I_b|\ge2\overline{p}+2|L|+4\ge12.
$$
This yields $$
|C|\ge12+4(\delta-4)=4\delta-4\ge3\delta+1,
$$
contradicting (23). Otherwise, $|L|=1$ and $M(I_1,I_2,...,I_s)\subseteq E(G)$. By Lemma 2, $|I_a|+|I_b|\ge2\overline{p}+6=10$, implying that $|C|\ge10+4(\delta-4)=4\delta-6$. Combining this with (23), we get $4\delta-6\le |C|\le3\delta-1$, i.e. $\delta\le5$. Since $\delta\ge5$, we have 
$$
\delta=5,  \  s=3,  \   |C|=3\delta-1=14,          \   |I_a|+|I_b|=10,
$$
$$
|I_i|=4   \    \mbox{for each}  \   i\in\{1,2,...,s\}\backslash \{a,b\}.
$$
Assume w.l.o.g. that $a=1$ and $b=2$. By Lemma 2, $|M(I_1,I_2)|=1$ and $M(I_1,I_3)=M(I_2,I_3)=\emptyset$, i.e. $|M(I_1,I_2,...,I_s)|=1$. Claim 3 is proved.    \quad      $\Delta$                    \\

Put $G^\prime=G\backslash M(I_1,I_2,...,I_s)$. As in Case 3.1.2, form a graph $G^\ast$ by adding in $G^\prime$ at most $2|M(I_1,I_2,...,I_s)|$ appropriate new edges such that $\delta(G^\ast)=\delta(G)$, $G^\ast\backslash \{\xi_1,\xi_2,...,\xi_s\}$ is disconnected and 
$$
q(G)\ge q(G^\ast)-|M(I_1,I_2,...,I_s)|.   
$$
Let $H_1,H_2,...,H_{s+1}$ be the connected components of $G^\ast\backslash \{\xi_1,\xi_2,...,\xi_s\}$ with $V(I_i^\ast)\subseteq V(H_i)$ $(i=1,2,...,s)$ and $V(H_{s+1})=\{x_1,x_2,x_3\}$. Using notation (5) for $G^\ast$, we have $h_i\ge3$ $(i=1,2,...,s+1)$. If $h_i\ge6$ for some $i\in \{1,2,...,s\}$ then 
$$
|C|\ge6+3(s-1)+|\{\xi_1,\xi_2,...,\xi_s\}|=4\delta-5\ge3\delta,
$$
contradicting (23). Let $3\le h_i\le 5$ $(i=1,2,...,s+1)$. By Lemma 3,
$$
q_i(G^\ast)\ge\frac{h_i(2\delta-h_i+1)}{2}\ge 3(\delta-1)   \  \  (i=1,2,...,s+1),
$$ 
implying that
$$
q(G^\ast)\ge\sum_{i=1}^{s+1}q_i(G^\ast)\ge3(s+1)(\delta-1)=3(\delta-1)^2.
$$
By Claim 3, 
$$
q\ge q(G^\ast)-1\ge3(\delta-1)^2-1>\frac{3(\delta-1)(\delta+2)}{2}.
$$
 
\textbf{Case 4.1.3}. $\overline{p}=3$.

Put $P=x_1x_4x_3x_2$. It follows that $|N_C(x_i)|\ge\delta-3\ge2$ $(i=1,2)$. By (1), $|C|\le3\delta+1-\overline{p}=3\delta-2$. Combining this with (16), we have
$$
3\delta-3\le |C|\le3\delta-2.           \eqno{(24)}
$$ 
If $N_C(x_1)\not= N_C(x_2)$ then by Lemma 1, 
$$
|C|\ge4\delta-2\overline{p}=4\delta-6\ge3\delta-1,
$$
contradicting (24). Let $N_C(x_1)= N_C(x_2)$. Clearly, $|I_i|\ge5$ $(i=1,2,...,s)$. If $s\ge \delta-2$ than $|C|\ge5s\ge5\delta-10>3\delta-1$, contradicting (24).  Hence $s\le \delta-3$. Observing also that $s=|N_C(x_1)|\ge\delta-3$, we get $s=\delta-3$. It follows that $x_1x_2\in E(G)$. By symmetric arguments, $N_C(x_i)=N_C(x_1)$ $(i=2,3,4)$.\\

\textbf{Claim 4}. $M(I_1,I_2,...,I_s)=\emptyset$.

\textbf{Proof}. Assume to the contrary, i.e. $M(I_1,I_2,...,I_s)\not=\emptyset$. It means that $M(I_a,I_b)\not=\emptyset$ for some distinct $a,b\in\{1,2,...,s\}$. By Lemma 2, 
$$
|I_a|+|I_b|\ge4\delta-2\overline{p}=4\delta-6,
$$
implying that $|C|\ge(4\delta-6)+5(s-2)=9\delta-31$. Combining this with (24), we get $9\delta-31\le|C|\le 3\delta-2$, a contradiction. Claim 4 is proved.   \quad     $\Delta$  \\

By Claim 4, $G\backslash \{\xi_1,\xi_2,...,\xi_s\}$ is disconnected. Let $H_1,H_2,...,H_{s+1}$ be the connected components of $G\backslash \{\xi_1,\xi_2,...,\xi_s\}$ with $V(I_i^\ast)\subseteq V(H_i)$ $(i=1,2,...,s)$ and $V(H_{s+1})=V(P)$. By notation (5), we have $h_i\ge4$ $(i=1,2,...,s+1)$. If $h_i\ge8$ for some $i\in\{1,2,...,s\}$ then 
$$
|C|\ge 8+4(\delta-4)+s=5\delta-11\ge3\delta-1,
$$
contradicting (24).  Let $4\le h_i\le7$ $(i=1,2,...,s+1)$. By Lemma 3, 
$$
q_i\ge\frac{h_i(2\delta-h_i+1)}{2}\ge2(2\delta-3)   \quad    (i=1,2,...,s+1).
$$
Hence 
$$
q\ge\sum_{i=1}^{s+1}q_i\ge2(s+1)(2\delta-3)
$$
$$
=2(\delta-2)(2\delta-3)\ge\frac{3(\delta-1)(\delta+2)}{2}.
$$

 \textbf{Case 4.1.4}. $\overline{p}=4$.

Put $P=x_1x_5x_4x_3x_2$. By (1), $|C|\le3\delta+1-\overline{p}=3\delta-3$, and by (16), $|C|\ge3\delta-3$. It follows that 
$$
|C|=3\delta-3, \ n=3\delta+2,  \   V(G)=V(C\cup P).            \eqno{(25)}
$$
If $\delta\le6$ then
$$
q\ge\frac{n\delta}{2}=\frac{(3\delta+2)\delta}{2}\ge\frac{3(\delta-1)(\delta+2)}{2}.
$$
Let $\delta\ge7$, implying that $|N_C(x_i)|>2$ $(i=1,2)$. If $N_C(x_1)\not= N_C(x_2)$ then by Lemma 1, 
$$
|C|\ge4\delta-2\overline{p}=4\delta-8\ge3\delta-2,
$$
contradicting (25).  Let $N_C(x_1)= N_C(x_2)$. Clearly, $|I_i|\ge6$ $(i=1,2,...,s)$. If $s\ge \delta-3$ then 
$$
|C|\ge(\overline{p}+2)s\ge6(\delta-3)\ge3\delta-2,
$$
contradicting (25). Let $s\le\delta-4$. On the other hand, $s\ge|N(x_1)|-\overline{p}\ge\delta-4$, implying that $s=\delta-4$. It follows that $x_1x_2\in E(G)$. By symmetric arguments, $N_C(x_i)=N_C(x_1)$ $(i=2,3,4,5)$. If $M(I_a,I_b)\not=\emptyset$ for some distinct elementary segments $I_a,I_b$, then by Lemma 2, 
$$
|I_a|+|I_b|\ge4\delta-2\overline{p}=4\delta-8.
$$
Hence
$$
|C|\ge4\delta-8+6(s-2)=10\delta-44\ge3\delta-2,
$$
contradicting (25). Otherwise, $M(I_1,I_2,...,I_s)=\emptyset$, i.e. $G\backslash \{\xi_1,\xi_2,...,\xi_s\}$ is disconnected. Let $H_1,H_2,...,H_{s+1}$ be the connected components of $G\backslash \{\xi_1,\xi_2,...,\xi_s\}$ with $V(I_i^\ast)\subseteq V(H_i)$ $(i=1,2,...,s)$ and $V(H_{s+1})=V(P)$. By notation (5), $h_i\ge5$ $(i=1,2,...,s+1)$. If $h_i\ge6$ for some $i\in \{1,2,...,s\}$ then
$$
|C|\ge6+5(s-1)+s=6\delta-23\ge3\delta-2,
$$
contradicting (25). So, $h_i=5$ $(i=1,2,...,s+1)$. By Lemma 3,
$$
q_i\ge\frac{h_i(2\delta-h_i+1)}{2}=5(\delta-2)   \   \   (i=1,2,...,s+1),
$$
implying that
$$
q\ge\sum_{i=1}^{s+1}q_i\ge5(s+1)(\delta-2)=5(\delta-3)(\delta-2)>\frac{3(\delta-1)(\delta-2)}{2}.
$$

\textbf{Case 4.2}. $\delta-2\le\overline{p}\le\delta+1$.

\textbf{Case 4.2.1}. $\overline{p}=\delta-2$.

It follows that $|N_C(x_i)|\ge2$ $(i=1,2)$. By (1), 
$$
|C|\le 3\delta+1-\overline{p}=2\delta+3.               \eqno{(26)}
$$
If $N_C(x_1)\not=N_C(x_2)$ then by Lemma 1, $|C|\ge4\delta-2\overline{p}=2\delta+4$, contradicting (26). Let $N_C(x_1)=N_C(x_2)$. If $s\ge3$ then $|C|\ge s(\overline{p}+2)\ge3\delta\ge2\delta+4$,  again contradicting (26). Let $s=2$. It follows that $x_1x_2\in E(G)$. By symmetric arguments, $N_C(y)=N_C(x_1)=\{\xi_1,\xi_2\}$ for each $y\in V(P)$. Clearly, $|I_i|\ge\overline{p}+2=\delta$ $(i=1,2)$.\\

\textbf{Case 4.2.1.1}. $M(I_1,I_2)=\emptyset$.

It follows that $G\backslash \{\xi_1,\xi_2\}$ is disconnected. Let $H_1,H_2,...,H_t$ be the connected components of $G\backslash \{\xi_1,\xi_2\}$ with $V(I_i^\ast)\subseteq V(H_i)$ $(i=1,2)$ and $V(P)\subset V(H_3)$. Since $G[V(P)]$ is hamiltonian, we have $V(H_3)=V(P)$. By notation (5),  $h_i\ge|I_i|-1\ge\delta-1$ $(i=1,2)$ and $h_3=\delta-1$. If $h_i\ge\delta+3$ for some $i\in \{1,2\}$ then 
$$
|C|\ge(\delta+3)+(\delta-1)+|\{\xi_1,\xi_2\}|=2\delta+4,
$$
contradicting (26). So, $\delta-1\le h_i\le \delta+2$ $(i=1,2,3)$.  By Lemma 3,
$$
q_i\ge\frac{h_i(2\delta-h_i+1)}{2}\ge \frac{(\delta-1)(\delta+2)}{2}   \quad     (i=1,2,3).
$$ 
Hence,  
$$
q\ge \sum_{i=1}^3q_i\ge \frac{3(\delta-1)(\delta+2)}{2}.
$$

\textbf{Case 4.2.1.2}. $M(I_1,I_2)\not=\emptyset$.

By the definition, there is an intermediate path $L$ between $I_1$ and $I_2$. If $|L|\ge2$ then by Lemma 2, 
$$
|C|=|I_1|+|I_2|\ge2\overline{p}+2|L|+4\ge2\delta+4,
$$
contradicting (26).  Otherwise, $M(I_1,I_2)\subseteq E(G)$. Further, if $|M(I_1,I_2)|\ge3$ then by Lemma 2, 
$$
|C|=|I_1|+|I_2|\ge2\overline{p}+8=2\delta+4,
$$
contradicting (26). Thus $|M(I_1,I_2)|\le2$.\\

\textbf{Case 4.2.1.2.1}. $|M(I_1,I_2)|=1$.

Put $G^\prime=G\backslash M(I_1,I_2)$. As in Case 3.1.2, form a graph $G^\ast$ by adding at most two new edges in $G^\prime$ such that $\delta(G^\ast)=\delta(G)$, $G^\ast\backslash \{\xi_1,\xi_2\}$ is disconnected and $q(G)\ge q(G^\ast)-1$. Let $H_1,H_2,...,H_t$ be the connected components of $G^\ast\backslash \{\xi_1,\xi_2\}$ with $V(I_i^\ast)\subseteq V(H_i)$ $(i=1,2)$ and $V(P)= V(H_3)$. Using notation (5) for $G^\ast$, as in Case 4.2.1.1, we have $\delta-1\le h_i\le \delta+2$ $(i=1,2,3)$. By Lemma 2, $|I_1|+|I_2|\ge2\overline{p}+6=2\delta+2$. It means that $\max_i |I_i|\ge \delta+1$, i.e. $\max_i h_i\ge\delta$. Assume w.l.o.g. that $h_1\ge\delta$. By Lemma 3,
$$
q_1(G^\ast)\ge\frac{h_1(2\delta-h_1+1)}{2}\ge \frac{\delta(\delta+1)}{2},   
$$ 

$$
q_i(G^\ast)\ge\frac{h_i(2\delta-h_i+1)}{2}\ge \frac{(\delta-1)(\delta+2)}{2}   \quad     (i=2,3),
$$ 
implying that 
$$
q(G^\ast)\ge\frac{\delta(\delta+1)}{2}+(\delta-1)(\delta+2).
$$
Hence,  
$$
q\ge q(G^\ast)-1\ge \frac{\delta(\delta+1)}{2}+(\delta-1)(\delta+2)-1\ge\frac{3(\delta-1)(\delta+2)}{2}.
$$

\textbf{Case 4.2.1.2.2}. $|M(I_1,I_2)|=2$.

By Lemma 2, 
$$
|C|= |I_1|+|I_2|\ge2\overline{p}+7=2\delta+3.
$$
By (26),  $|C|=2\delta+3$ and $V(G)=V(P\cup C)$.
 Put $G^\prime=G\backslash M(I_1,I_2)$. As in Case 3.1.2, form a graph $G^\ast$ by adding at most four new edges in $G^\prime$ such that 
$\delta(G^\ast)=\delta(G)$, $G^\ast\backslash \{\xi_1,\xi_2\}$ is disconnected and $q(G)\ge q(G^\ast)-2$. Let $H_1,H_2,H_3$ 
be the connected components of $G^\ast\backslash \{\xi_1,\xi_2\}$ with $V(H_i)=V(I_i^\ast)$ $(i=1,2)$ and $V(H_3)=V(P)$.
 Using notation (5) for $G^\ast$, we have as in Case 4.2.1.1, $\delta-1\le h_i \le \delta+2$ $(i=1,2,3)$. 
Since $|I_i|\ge\delta$ and $|C|=|I_1|+|I_2|=2\delta+3$, we can assume w.l.o.g. that either $|I_1|=\delta+2$, $|I_2|=\delta+1$ or $|I_1|=\delta+3$, $|I_2|=\delta$.\\

\textbf{Case 4.2.1.2.2.1}. $|I_1|=\delta+2$, $|I_2|=\delta+1$.

It follows that $h_1=\delta+1$, $h_2=\delta$ and $h_3=\delta-1$. By Lemma 3,

$$
q_i(G^\ast)\ge\frac{h_i(2\delta-h_i+1)}{2}= \frac{\delta(\delta+1)}{2}   \quad     (i=1,2),
$$ 
$$
q_3(G^\ast)\ge\frac{h_3(2\delta-h_3+1)}{2}=\frac{(\delta-1)(\delta+2)}{2}.
$$
Hence
$$
q\ge \sum_{i=1}^3q_i(G^\ast)-2= \frac{3(\delta-1)(\delta+2)}{2}.
$$

\textbf{Case 4.2.1.2.2.2}. $|I_1|=\delta+3$, $|I_2|=\delta$.

Let $M(I_1,I_2)=\{e_1,e_2\}$, where 
$$
e_1=y_1z_1,   \    e_2=y_2z_2,  \  \{y_1,y_2\}\subseteq V(I_1^\ast), \  \{z_1,z_2\}\subseteq V(I_2^\ast).
$$
If $y_1\not=y_2$ and $z_1\not=z_2$ then by Lemma 2, 
$$
|I_1|+|I_2|\ge2\overline{p}+8=2(\delta-2)+8=2\delta+4,
$$
contradicting (26). Let either $y_1\not=y_2$ and $z_1=z_2$ or $y_1=y_2$ and $z_1\not=z_2$.\\

\textbf{Case 4.2.1.2.2.2.1}. $y_1\not=y_2$ and $z_1=z_2$.

Assume w.l.o.g. that $y_1,y_2$ occur on $I_1$ in this order. If $y_2=y_1^+$ then
$$
|C|\ge |\xi_1\overrightarrow{C}y_1z_1y_2\overrightarrow{C}\xi_2x_2\overleftarrow{P}x_1\xi_1|=2\delta+4,
$$
contradicting (26). Let $y_2\not=y_1^+$, i.e. $|y_1\overrightarrow{C}y_2|\ge2$. Put 
$$
C^\prime=\xi_1\overrightarrow{C}y_2z_1\overleftarrow{C}\xi_2x_2\overleftarrow{P}x_1\xi_1,
$$
$$
C^{\prime\prime}=\xi_1\overleftarrow{C}z_1y_1\overrightarrow{C}\xi_2x_2\overleftarrow{P}x_1\xi_1.
$$
Clearly,
$$
|C|\ge|C^\prime|=|\xi_1\overrightarrow{C}y_1|+|y_1\overrightarrow{C}y_2|+1+|\xi_2\overrightarrow{C}z_1|+\overline{p}+2,
$$
$$
|C|\ge|C^{\prime\prime}|=|\xi_1\overleftarrow{C}z_1|+|y_1\overrightarrow{C}y_2|+|y_2\overrightarrow{C}\xi_2|+1+\overline{p}+2.
$$
By summing, we get 
$$
2|C|\ge(|\xi_1\overrightarrow{C}y_1|+|y_1\overrightarrow{C}y_2|+|y_2\overrightarrow{C}\xi_2|+|\xi_2\overrightarrow{C}z_1|+|z_1\overrightarrow{C}\xi_1|)+|y_1\overrightarrow{C}y_2|+2+2\delta
$$
$$
=|C|+|y_1\overrightarrow{C}y_2|+2\delta+2\ge|C|+2\delta+4.
$$
Hence $|C|\ge2\delta+4$, contradicting (26).\\

\textbf{Case 4.2.1.2.2.2.2}. $y_1=y_2$ and $z_1\not=z_2$. 

Assume w.l.o.g. that $z_2, z_1$ occur on $I_2$ in this order.\\

\textbf{Case 4.2.1.2.2.2.2.1}. $\delta\ge6$.

If $|\xi_1\overrightarrow{C}y_1|\ge\delta-1$ and $|y_1\overrightarrow{C}\xi_2|\ge\delta-1$ then $|I_1|\ge2\delta-2\ge\delta+4$, contradicting the hypothesis. Thus, we can assume w.l.o.g. that $|\xi_1\overrightarrow{C}y_1|\le\delta-2$. If $y_1^-=\xi_1$ then 
$$
|\xi_1\overleftarrow{C}z_2y_1\overrightarrow{C}\xi_2x_2\overleftarrow{P}x_1\xi_1|\ge2\delta+5,
$$
contradicting (26). Let $y_1^-\not=\xi_1$, that is $y_1^-\in V(I_1^\ast)$. Since $|M(I_1,I_2)|=2$, we have $N(y_1^-)\subset V(I_1)$. If $N(y_1^-)\cap V(y_1^+\overrightarrow{C}\xi_2^-)=\emptyset$ then $|N(y_1^-)|\le\delta-1$, a contradiction. Otherwise, $y_1^-w\in E(G)$ for some $w\in V(y_1^+\overrightarrow{C}\xi_2^-)$. Put 
$$
R=\xi_1\overrightarrow{C}y_1^-w\overleftarrow{C}y_1
$$
$$
C^\prime=\xi_1\overrightarrow{R}y_1z_1\overleftarrow{C}\xi_2x_2\overleftarrow{P}x_1\xi_1,
$$
$$
C^{\prime\prime}=\xi_1\overleftarrow{C}z_2y_1\overrightarrow{C}\xi_2x_2\overleftarrow{P}x_1\xi_1.
$$
Clearly, $|R|\ge|\xi_1\overrightarrow{C}y_1|+1$ and
$$
|C|\ge|C^\prime|=|R|+1+|z_1\overleftarrow{C}\xi_2|+(\overline{p}+2)\ge|\xi_1\overrightarrow{C}y_1|+|z_1\overleftarrow{C}\xi_2|+\delta+2,
$$
$$
|C|\ge|C^{\prime\prime}|=|\xi_1\overleftarrow{C}z_1|+2+|y_1\overrightarrow{C}\xi_2|+(\overline{p}+2).
$$
By summing, we get
$$
2|C|\ge(|\xi_1\overrightarrow{C}y_1|+|y_1\overrightarrow{C}\xi_2|+|\xi_2\overrightarrow{C}z_1|+|z_1\overrightarrow{C}\xi_1|)+2\delta+4=|C|+2\delta+4.
$$
Hence $|C|\ge2\delta+4$, contradicting (26).\\

\textbf{Case 4.2.1.2.2.2.2.2}. $\delta=5$.

It follows that 
$$
|I_1|=\delta+3=8,  \   |I_2|=\delta=5,           \        |C|=2\delta+3=13.
$$
If either $|\xi_1\overrightarrow{C}y_1|\le\delta-2=3$ or $|y_1\overrightarrow{C}\xi_2|\le\delta-2=3$ then we can argue as in Case 4.2.1.2.2.2.2.1. Otherwise, $|\xi_1\overrightarrow{C}y_1|=|y_1\overrightarrow{C}\xi_2|=4$. If $|z_1\overleftarrow{C}\xi_2|\ge4$ then
$$
|\xi_1\overrightarrow{C}y_1z_1\overleftarrow{C}\xi_2x_2\overleftarrow{P}x_1\xi_1|\ge14>|C|,
$$
a contradiction. Let $|z_1\overleftarrow{C}\xi_2|\le3$. Analogously, $|\xi_1\overleftarrow{C}z_2|\le3$, implying that 
$I_2=\xi_1z_1^+z_1z_2\xi_2^+\xi_2$.
 If $z_1^+z_2\in E(G)$ then
$$
|\xi_1\overrightarrow{C}y_1z_1z_1^+z_2\overleftarrow{C}\xi_2x_2\overleftarrow{P}x_1\xi_1|=14>|C|,
$$
a contradiction.  So, $N(z_1^+)\subseteq \{\xi_1,\xi_2,z_1,\xi_2^+\}$, again a contradiction, since $|N(z_1^+)|\ge\delta=5$.\\

\textbf{Case 4.2.2}. $\overline{p}=\delta-1$.

By (1),
$$
|C|\le3\delta+1-\overline{p}=2\delta+2.                \eqno{(27)}
$$
It follows that $|N_C(x_i)|\ge1$ $(i=1,2)$.\\

\textbf{Case 4.2.2.1}. $|N_C(x_i)|\ge2$ $(i=1,2)$.

If $N_C(x_1)\not=N_C(x_2)$ then by Lemma 1, 
$|C|\ge2\overline{p}+8=2\delta+6$, contradicting (27). 
Let $N_C(x_1)=N_C(x_2)$. If $s\ge3$ then 
$$
|C|\ge s(\overline{p}+2)\ge3(\delta+1)>2\delta+2,
$$
contradicting (27). Let $s=2$. It follows that
$$
|C|=2\delta+2,   \       |I_1|=|I_2|=\delta+1,        \     V(G)=V(C\cup P).
$$
Assume that $yz\in E(G)$ for some  $y\in V(P)$ and $z\in V(C)\backslash \{\xi_1,\xi_2\}$. Assume w.l.o.g. that $z\in V(I_1^\ast)$. Since $\overline{p}=\delta-1\ge4$, we can assume w.l.o.g. that $|x_1\overrightarrow{P}y|\ge2$. If $x_2w\in E(G)$ for some $w\in \{y^-,y^{-2}\}$ then
$$
|\xi_1\overrightarrow{C}z|\ge|\xi_1x_1\overrightarrow{P}wx_2\overleftarrow{P}yz|\ge\delta.
$$
Observing also that $|z\overrightarrow{C}\xi_2|\ge2$, we have $|I_1|\ge\delta+2$, a contradiction. Otherwise, $d(x_2)\le\delta-1$, a contradiction. So, $N_C(y)\subseteq \{\xi_1,\xi_2\}$ for each $y\in V(P)$. On the other hand, by Lemma 2, $M(I_1,I_2)=\emptyset$ and hence $G\backslash \{\xi_1,\xi_2\}$ is disconnected. Let $H_1,H_2,H_3$ be the connected components of $G\backslash \{\xi_1,\xi_2\}$ with $V(H_i)=V(I_i^\ast)$ $(i=1,2)$ and $V(H_3)=V(P)$. By notation (5),  $h_i=\delta$ $(i=1,2,3)$. By Lemma 3, 
$$
q_i\ge\frac{h_i(2\delta-h_i+1)}{2}=\frac{\delta(\delta+1)}{2}   \quad       (i=1,2,3),
$$
implying that
$$
q\ge\sum_{i=1}^3q_i\ge\frac{3(\delta^2+\delta)}{2}>\frac{3(\delta-1)(\delta+2)}{2}.
$$

\textbf{Case 4.2.2.2}. Either $|N_C(x_1)|=1$ or $|N_C(x_2)|=1$.

Assume w.l.o.g. that $|N_C(x_1)|=1$. It follows that $V(P)\backslash \{x_1\}\subset N(x_1)$. Put $N_C(x_1)=\{y_1\}$.\\

 \textbf{Case 4.2.2.2.1}. $N_C(x_2)\not=N_C(x_1)$.

It follows that $x_2y_2\in E(G)$ for some $y_2\in V(C)\backslash \{y_1\}$. Clearly, $|y_1\overrightarrow{C}y_2|\ge\delta+1$ and $|y_2\overrightarrow{C}y_1|\ge\delta+1$. Hence
$$
|y_1\overrightarrow{C}y_2|=|y_2\overrightarrow{C}y_1|=\delta+1,   \      |C|=2\delta+2,        \         V(G)=V(C\cup P).                              \eqno{(28)}
$$ 
If $s\ge3$ then there are at least two elementary segments on $C$ of length at least $\delta+1$. It means that $|C|>2\delta+2$, contradicting (28). Let $s=2$, i.e. $N_C(x_1)\cup N_C(x_2)=\{y_1,y_2\}=\{\xi_1,\xi_2\}$. Assume that $zw\in E(G)$ for some $z\in V(P)$ and $w\in V(C)\backslash \{y_1,y_2\}$, and assume w.l.o.g. that $w\in y_1\overrightarrow{C}y_2$. Since $V(P)\backslash \{x_1\}\subset N(x_1)$, we have $x_1z^+\in E(G)$. Further, since $C$ is extreme, 
$$
|w\overrightarrow{C}y_2|\ge|wz\overleftarrow{P}x_1z^+\overrightarrow{P}x_2y_2|\ge\delta+1.
$$
Hence, $|C|>2\delta+2$, contradicting (28). Thus, $N(z)\subseteq \{y_1,y\}$ for each $z\in V(P)$. On the other hand, by Lemma 2, $M(I_1,I_2)=\emptyset$. Further, we can argue as in Case 4.2.2.1.\\

 \textbf{Case 4.2.2.2.2}. $N_C(x_2)=N_C(x_1)=\{y_1\}$.

It follows that 
$$
N(x_i)=(V(P)\backslash \{x_i\})\cup \{y_1\}  \   \   (i=1.2).
$$
Since $\kappa\ge2$, there is a path $R=z\overrightarrow{R}w$ such that $z\in V(P)$ and $w\in V(C)\backslash\{y_1\}$. Since $N_C(x_1)=N_C(x_2)=\{y_1\}$, we have $z\not\in \{x_1,x_2\}$. Then
$$
|y_1x_1\overrightarrow{P}z^-x_2\overleftarrow{P}zw|=\delta+1,
$$
and we can argue as in Case 4.2.2.2.1.\\

\textbf{Case 4.2.3}. $\overline{p}=\delta$.

By (1), $|C|\le 3\delta+1-\overline{p}=2\delta+1$.  If $|Q|\ge\delta+1$ then by (2), $|C|\ge2|Q|\ge2\delta+2$, a contradiction. Let 
$$
|Q|\le\delta.                              \eqno{(29)}
$$

\textbf{Case 4.2.3.1}. $x_1x_2\not\in E(G)$

It follows that $|N_C(x_i)|\ge1$ $(i=1,2)$. If  $|N_C(x_i)|\ge2$ for some  $i\in\{1,2\}$ then clearly 
$|Q|\ge\overline{p}+2=\delta+2$, contradicting (29). Let $|N_C(x_1)|=|N_C(x_2)|=1$.
 Further, if $N_C(x_1)\not=N_C(x_2)$
 then again $|Q|\ge\delta+2$,  contradicting (29).
 Let $N_C(x_1)=N_C(x_2)=\{z_1\}$ for some
 $z_1\in V(C)$. Since $\kappa\ge2$, 
there is a path $L=yz_2$
 connecting $P$ and $C$ such that 
$y\in V(P)$ and $z_2\in V(C)\backslash \{z_1\}$. 
Clearly, $y\not\in \{x_1,x_2\}$. If  $x_2y^-\in E(G)$ then
$$
|Q|\ge|z_1x_1\overrightarrow{P}y^-x_2\overleftarrow{P}yz_2|=\delta+2,
$$
contradicting (29). Let  $x_2y^-\not\in E(G)$. Further, if $y^-\not=x_1$ then recalling that $x_2x_1\not\in E(G)$, we conclude that $|N_C(x_2)|\ge2$, a contradiction. Otherwise, $y^-=x_1$ and $|Q|\ge|z_1x_2\overleftarrow{P}yz_2|=\delta+1$, contradicting (29).\\

\textbf{Case 4.2.3.2}. $x_1x_2\in E(G)$.

Put $C^\prime=x_1\overrightarrow{P}x_2x_1$. Since $\kappa\ge2$, there are two disjoint paths $L_1,L_2$ connecting $C^\prime$ and $C$. Further, since $P$ is extreme, $|L_1|=|L_2|=1$. Let $L_1=y_1z_1$ and $L_2=y_2z_2$, where, $y_1,y_2\in V(C^\prime)$ and $z_1,z_2\in V(C)$. Since $C^\prime$ is a Hamilton cycle in $G[V(P)]$ and $|C^\prime|\ge\delta+1\ge6$, we can assume that $P$ is chosen such that $x_1=y_1$ and $|x_1\overrightarrow{P}y_2|\ge3$. If  $x_2v\in E(G)$ for some $v\in \{y_2^{-1},y_2^{-2}\}$ then
$$
|Q|\ge|z_1x_1\overrightarrow{P}vx_2\overleftarrow{P}y_2z_2|\ge\delta+1,
$$
contradicting (29). Otherwise, $|N_C(x_2)|\ge2$, implying that $x_2z_3\in E(G)$ for some $z_3\in V(C)\backslash \{z_1\}$. Then 
$$
|Q|\ge|z_1x_1\overrightarrow{P}x_2z_3|\ge\delta+2,
$$
again contradicting (29).\\

\textbf{Case 4.2.4}. $\overline{p}=\delta+1$.

By (1), $|C|\le 3\delta+1-\overline{p}=2\delta$. Recalling (14), we get $|C|=2\delta$ and $V(G)=V(C\cup P)$.  If $|Q|\ge\delta+1$ then by (2), $|C|\ge2|Q|\ge2\delta+2$, a contradiction. Let 
$$
|Q|\le\delta.                                      \eqno{(30)}
$$

\textbf{Case 4.2.4.1}. $x_1x_2\in E(G)$.

Put $C^\prime=x_1\overrightarrow{P}x_2x_1$. Since $\kappa\ge2$, there are two disjoint edges $z_1w_1$ and $z_2w_2$ connecting $C^\prime$ and $C$ such that $z_1,z_2\in V(C^\prime)$ and $w_1,w_2\in V(C)$. Since $C^\prime$ is a Hamilton cycle in $G[V(P)]$ and $|C^\prime|\ge\delta+2\ge7$, we can assume w.l.o.g. that $P$ is chosen such that $z_1=x_1$ and $|x_1\overrightarrow{P}z_2|\ge4$. If $x_2v\in E(G)$ for some $v\in \{z_2^{-1},z_2^{-2},z_2^{-3}\}$ then
$$
|Q|\ge|w_1x_1\overrightarrow{P}vx_2\overleftarrow{P}z_2w_2|\ge\delta+1,
$$
contradicting (30). Now let  $x_2v\not\in E(G)$ for each $v\in \{z_2^{-1},z_2^{-2},z_2^{-3}\}$. It follows that $|N_C(x_2)|\ge2$, i.e. $x_2w_3\in E(G)$ for some $w_3\in VC)\backslash \{w_1\}$. But then $|Q|\ge|w_1x_1\overrightarrow{P}x_2w_3|=\delta+3$, contradicting (30).\\

\textbf{Case 4.2.4.2}. $x_1x_2\not\in E(G)$.

If $d_P(x_1)+d_P(x_2)\ge|V(P)|=\overline{p}+2$ then by Theorem F, $G[V(P)]$ is hamiltonian and we can argue as in Case 4.2.4.1. Otherwise,  
$$
d_C(x_1)+d_C(x_2)\ge \delta-1\ge4.                  \eqno{(31)}
$$
Assume w.l.o.g. that $d_C(x_1)\ge d_C(x_2)$. \\

\textbf{Case 4.2.4.2.1}. $d_C(x_2)=0$.

It follows that $N(x_2)=V(P)\backslash \{x_2\}$. By (31),  $d_C(x_1)\ge4$. Put $C^\prime=x_1^+\overrightarrow{P}x_2x_1^+$. Since $\kappa\ge2$, there is a path $L=z\overrightarrow{L}w$ connecting $C^\prime$ and $C$ such that $z\in V(C^\prime)\backslash \{x_1^+\}$ and $w\in V(C)$. If $x_1\in V(L)$, i.e. $x_1z\in E(G)$, then $x_1\overrightarrow{P}z^-x_2\overleftarrow{P}zx_1$ is a Hamilton cycle in $G[V(P)]$ and we can argue as in Case 4.2.4.1. Let $x_1\not\in V(L)$. Since $V(G)=V(C\cup P)$, we have $L=zw$. Further, since $d_C(x_1)\ge4$, we have $x_1w_1\in E(G)$ for some $w_1\in V(C)\backslash \{w\}$. Hence, 
$$
|Q|\ge|w_1x_1\overrightarrow{P}z^-x_2\overleftarrow{P}zw|=\delta+3,
$$
contradicting (30).\\

\textbf{Case 4.2.4.2.2}. $d_C(x_2)=1$.

Let $N_C(x_2)=\{w_1\}$. By (31), $d_C(x_1)\ge3$, i.e. $x_1w\in E(G)$ for some $w\in V(C)\backslash \{w_1\}$. Hence 
$$
|Q|\ge|wx_1\overrightarrow{P}x_2w_1|=\delta+3,
$$
contradicting (30).\\

\textbf{Case 4.2.4.2.3}. $d_C(x_2)\ge2$.

Since $d_C(x_1)\ge d_C(x_2)$, we have $d_C(x_1)\ge2$. Hence $|Q|\ge\overline{p}+2=\delta+3$, contradicting (30).         \quad         \quad       \rule{7pt}{6pt}            \\

\noindent\textbf{Proof of Theorem 1}. Let $G$ be a 2-connected graph, $C$ a longest cycle in $G$ and $P=x_1\overrightarrow{P}x_2$ a longest path in $G\backslash C$ of length $\overline{p}$. If $\overline{p}=0$ then $C$ is a dominating cycle and we are done. Let $\overline{p}\ge1$.\\

\textbf{Case 1}. $\delta=2$ and $q\le8$.

Since $\kappa\ge2$ and  $\overline{p}\ge1$, there exist a path $Q=\xi\overrightarrow{Q}\eta$ such that $|Q|\ge3$ and $V(Q)\cap V(C)=\{\xi,\eta\}$. Further, since $C$ is extreme, we have $|C|=|y\overrightarrow{C}z|+|z\overrightarrow{C}y|\ge2|Q|\ge6$ and therefore, $q\ge|C|+|Q|\ge9$, contradicting the hypothesis.\\

\textbf{Case 2}. $\delta\ge3$ and $q\le(3(\delta-1)(\delta+2)-1)/2$.

Since 
$$
q=\frac{1}{2}\sum_{u\in V(G)}d(u)\ge \frac{\delta n}{2},
$$
we have $\delta n/2 \le (3(\delta-1)(\delta+2)-1)/2$, which is equivalent to 
$$
\delta\ge \frac{n-2}{3}-\frac{1}{3}+\frac{7}{3\delta}.
$$
If $n=3t$ for some integer $t$, then 
$$
\delta\ge \frac{3t-2}{3}-\frac{1}{3}+\frac{7}{3\delta}=t-1+\frac{7}{3\delta},
$$
implying that $\delta\ge t=n/3>(n-2)/3$. If $n=3t+1$ for some integer $t$, then 
$$
\delta\ge \frac{3t-1}{3}-\frac{1}{3}+\frac{7}{3\delta}=t-\frac{2}{3}+\frac{7}{3\delta},
$$
implying that $\delta\ge t=(n-1)/3>(n-2)/3$. Finally, if $n=3t+2$ for some integer $t$, then 
$$
\delta\ge \frac{3t}{3}-\frac{1}{3}+\frac{7}{3\delta}=t-\frac{1}{3}+\frac{7}{3\delta},
$$
implying that $\delta\ge t=(n-2)/3$. So, $\delta\ge(n-2)/3$, in any case. By Lemma 4, each longest cycle in $G$ is a dominating cycle.          \quad           \rule{7pt}{6pt}

\noindent Institute for Informatics and Automation Problems\\ National Academy of Sciences\\
P. Sevak 1, Yerevan 0014, Armenia\\ E-mail: zhora@ipia.sci.am


\begin{thebibliography}{10}

\bibitem{[1]}	J.A. Bondy and U.S.R. Murty, Graph Theory with Applications. Macmillan, London and Elsevier, New York (1976).

\bibitem{[2]}  G. A. Dirac, Some theorems on abstract graphs, Proc. London, Math. Soc., 2 (1952) 69--81.

\bibitem{[3]}   P. Erd\"{o}s and T. Gallai, On maximal paths and circuits of graphs, Acta Math. Acad. Sci. Hungar. 10 (1959) 337-356.

\bibitem{[4]}  Zh. G. Nikoghosyan, Path-Extensions and Long Cycles in Graphs, Transactions of the Institute for Informatics and Automation Problems of the
NAS of RA and of the Yerevan State University, Mathematical Problems of Computer Science 19 (1998), 25-31.

\bibitem{[5]}	Zh.G. Nikoghosyan, A Size Bound for Hamilton Cycles, preprint, arXiv:1107.2201v1.

\bibitem{[6]}   O. Ore, Note on Hamilton circuits, Amer. Math. Monthly 67 (1960) 55.

\end{thebibliography}
\end{document}